\documentclass[a4paper,12pt]{amsart}
\usepackage[margin=2.5cm]{geometry}
\usepackage{amsmath,amsthm,amssymb,amsfonts,mathrsfs,mathtools,bm,tikz,amscd}
\usepackage[T1]{fontenc}
\usepackage{enumerate}
\usepackage{lmodern}
\usepackage[utf8]{inputenc}
\usepackage[numbers]{natbib}
\usepackage{hyperref}
\usepackage{url}
\usepackage{xcolor}
\allowdisplaybreaks
\numberwithin{equation}{section}
\usetikzlibrary{arrows.meta}

\theoremstyle{definition}
\newtheorem{dfn}{Definition}[section]

\newtheorem{rmk}[dfn]{Remark}

\theoremstyle{plain} 
\newtheorem{thm}[dfn]{Theorem}
\newtheorem{prop}[dfn]{Proposition}

\newtheorem{lem}[dfn]{Lemma}

\newcommand{\N}{\mathbb{N}}

\newcommand{\R}{\mathbb{R}}

\usepackage{enumitem}

\title[Kelvin-invariant positive solutions for critical elliptic systems]{Kelvin-invariant positive solutions for critical elliptic systems with repulsive interaction}

\author{Yohei Sato}
\address{Department of Mathematics, Saitama University, Shimo-Okubo 255, Sakura-ku Saitama-shi, 338-8570, JAPAN}
\email{ ysato@rimath.saitama-u.ac.jp}

\author{Reo Suzuki}
\address{Department of Mathematics, Saitama University, Shimo-Okubo 255, Sakura-ku Saitama-shi, 338-8570, JAPAN}
\email{ suzuki.r.990@ms.saitama-u.ac.jp}

\begin{document}

\begin{abstract}
In this paper, we consider the following critical nonlinear elliptic system with repulsive interaction:
\begin{equation*}
\tag{$\mathcal{P}_\beta$}
\left\{
\begin{aligned}
	- & \Delta u_1 = \mu_1 u_1^3 - \beta(x) u_1 u_2^2 \quad \text{in } \R^4,\\
	- & \Delta u_2 = \mu_2 u_2^3 - \beta(x) u_1^2 u_2 \quad \text{in } \R^4,\\
	& u_1, u_2 \in \mathcal{D}^{1,2}(\R^4),
\end{aligned}
\right.
\end{equation*}
where $\mu_1, \mu_2 > 0$ and $\beta(x) \in C(\R^4,\R)$ is a nonnegative function.
Assuming that $\beta$ is invariant under the map $x \mapsto -\frac{x}{\vert{}x\vert{}^2}$,
we establish the existence of positive solutions to $(\mathcal{P}_\beta)$ that are invariant under the Kelvin transform subject to suitable assumptions.
This result provides a critical counterpart to the work of Sato  [Y. Sato, \textit{Nonlinear Anal.} \textbf{234} (2023), 113302], 
who established the existence of an even ground state for the corresponding subcritical problem.
Furthermore, we investigate the asymptotic behavior of the obtained solutions as $\|\beta\|_\infty \to 0$.

\medskip

\noindent
MSC2020: 
35J57, 
35J50, 
35J91, 

\noindent
Keywords: 
coupled nonlinear Schr\"odinger equations; 
Schr\"odinger systems; 
repulsive interaction;
variational methods; 
Kelvin transform; 
critical exponent;

\end{abstract}

\maketitle

\section{Introduction}
In this paper, we consider the following critical nonlinear 
elliptic system with repulsive interaction:
\begin{equation}
\label{P}
\tag{$\mathcal{P}_\beta$}
  \left\{ 
\begin{aligned}    
- & \Delta u_1 = \mu_1 u_1^3 - \beta(x) u_1 u_2^2 \quad \text{in } \R^4,\\
- & \Delta u_2 = \mu_2 u_2^3 - \beta(x) u_1^2 u_2 \quad \text{in } \R^4,\\
& u_1, u_2 \in \mathcal{D}^{1,2}(\R^4).
\end{aligned} 
  \right.
\end{equation}
Here $\mathcal{D}^{1,2}(\R^4) = \overline{C_0^\infty(\R^4)}^{\|\nabla \cdot\|_{L^2(\R^4)}}$. 
We assume that $\mu_1, \mu_2 > 0$ and $\beta(x) \in C(\R^4,\R)$ satisfies the following conditions. 

\begin{itemize}
\item[($\beta$1)]  $\beta(x) \ge 0 \quad (x \in \R^4)$.
\item[($\beta$2)] $\displaystyle \beta\left(-\frac{x}{|x|^2}\right) = \beta(x)  \quad (x \in \R^4\setminus \{0\}) $.
\end{itemize}
Under assumption ($\beta$2), $\beta(x)$ satisfies 
\begin{equation*}
	\lim_{|x| \to \infty} \beta(x) = \beta(0).
\end{equation*}

For $1 \le p < \infty$, we use the following notation.
\begin{align*}
&\|u\|_p = \left(\int_{\R^4} |u|^p \, dx\right)^{\frac{1}{p}} , \quad 
\|u\|_{\infty} = \operatorname*{ess\,sup}_{x \in \R^4} |u(x)| ,\\
&(u,v)_{\mathcal{D}^{1,2}} = \int_{\R^4} \nabla u \cdot \nabla v \, dx , \quad
\|u\|_{\mathcal{D}^{1,2}} = (u,u)_{\mathcal{D}^{1,2}}^{\frac{1}{2}} .
\end{align*}
Also, let $\mathcal{D} := \mathcal{D}^{1,2}(\R^4)^2$ be endowed with the inner product
\begin{equation*}
((u_1,u_2),(v_1,v_2))_{\mathcal{D}} = (u_1,v_1)_{\mathcal{D}^{1,2}} + (u_2,v_2)_{\mathcal{D}^{1,2}},
\end{equation*}
and the associated norm $\|\cdot\|_{\mathcal{D}}$. In what 
follows, whenever we write an element $\bm{u} \in \mathcal{D}$, we mean $\bm{u} = (u_1,u_2)$. 
The associated energy functional \(I:\mathcal D\to\mathbb R\) is defined by 
\begin{equation*}
	 I(\bm{u}) = \frac{1}{2}\|\bm{u}\|_{\mathcal{D}}^2 - \frac{\mu_1}{4}\|u_1\|_4^4 - \frac{\mu_2}{4}\|u_2\|_4^4
	 + \frac{1}{2} B(\bm{u}), 
	\text{ where } B(\bm{u}) = \int_{\R^4} \beta(x) u_1^2 u_2^2 \, dx.
\end{equation*}
The system \eqref{P} is called repulsive if $\beta(x)$ is nonnegative and $\beta(x) \not\equiv 0$.

We first recall the following subcritical counterpart of $(\mathcal{P}_\beta)$:
\begin{equation}
	\label{1.1}
	\left\{ \begin{aligned}    
	- & \Delta u_1 + \lambda_1 u_1 = \mu_1 u_1^3 - \beta(x) u_1 u_2^2 \quad \text{in } \R^N,\\
	- & \Delta u_2 + \lambda_2 u_2 = \mu_2 u_2^3 - \beta(x) u_1^2 u_2 \quad \text{in } \R^N,\\
	& u_1, u_2 \in H^1(\R^N),
	\end{aligned} \right.
\end{equation}
where $N \le 3$, $\lambda_1, \lambda_2, \mu_1, \mu_2 > 0$, and $\beta(x) \in C(\R^N,\R)$ is nonnegative and $\beta(x) \not\equiv 0$. 
For constant coefficients, Lin and Wei \cite{Lin--Wei} showed that the system \eqref{1.1} admits no ground states with both components nontrivial in the repulsive case.
More precisely, let $\widetilde{I}:H^1(\R^N)^2 \to \R$ be the functional associated with \eqref{1.1}, and define 
\begin{equation*}
	\mathcal{M}=\left\{(u_1,u_2)\in H^1(\R^N)^2 \mid \widetilde{I}'(u_1,u_2)(u_1,0) = \widetilde{I}'(u_1,u_2)(0,u_2)=0, u_1\not=0, u_2\not=0 \right\}.
\end{equation*}
Then the infimum 
\begin{equation*}
	\inf_{(u_1,u_2)\in \mathcal{M}}\widetilde{I}(u_1,u_2)
\end{equation*}
is not attained. 
This lack of attainment is related to the noncompactness of the embedding $H^1(\R^N)\hookrightarrow L^4(\R^N)$ caused by translations in $\R^N$.
On the other hand,  Sirakov \cite{Sirakov} proved the existence of a ground state in the class of radially symmetric functions for $N=2,3$. 
That is, the infimum of $\widetilde{I}$ over $\mathcal{M}\cap(H^1_{\rm rad}(\R^N))^2$ is attained. 
A key ingredient is the compactness of the embedding $H_{\mathrm{rad}}^1(\R^N)\hookrightarrow L^4(\R^N)$ for $N=2,3$.
More recently, in \cite{Sato23}, we established the existence of a ground state in the class of even functions 
under suitable assumptions on the even function $\beta$. 
Unlike the radial case, the embedding $H_{\rm even}^1(\R^N)\hookrightarrow L^4(\R^N)$ is not compact. 
Nevertheless, we showed that the infimum of $\widetilde{I}$ over $\mathcal{M}\cap(H^1_{\rm even}(\R^N))^2$ is attained.
Thus, radial symmetry can be relaxed to even symmetry.

The aim of this paper is to establish a critical counterpart of our previous result \cite{Sato23}.
However, unlike the subcritical case, the lack of compactness of the embedding $\mathcal{D}^{1,2}(\R^4)\hookrightarrow L^4(\R^4)$ arises not only from translations but also from dilations.
In particular, even symmetry alone does not prevent the loss of compactness caused by dilations.
In \cite{Sato--Suzuki}, we suggested that invariance under the Kelvin transform can play a role analogous to even symmetry in the subcritical setting. 
More precisely, using a profile decomposition theorem, we established compactness of suitable Palais--Smale sequences below certain energy thresholds.
Motivated by this observation, we impose assumption $(\beta2)$ and investigate the existence of least energy solutions to $(\mathcal{P}_\beta)$ in the class of invariant under the Kelvin transform.

Various existence results for critical elliptic systems with repulsive interaction have been obtained by imposing suitable symmetry conditions.
Chen, Medina, and Pistoia \cite{Chen--Medina--Pistoia} constructed nonradial solutions for a critical elliptic system with a small constant repulsive interaction by means of a Lyapunov--Schmidt reduction.
Clapp and Pistoia \cite{Clapp--Pistoia} studied critical competitive elliptic systems in $\R^N$ with constant coupling coefficients (also see \cite{Clapp--Szulkin}). 
They imposed invariance under a group $\Gamma\subset O(N+1)$ acting conformally on $\R^N$ through stereographic projection. 
By choosing $\Gamma$ so that every orbit on $\mathbb{S}^N$ has positive dimension, they obtained compactness in the corresponding invariant function space and established the existence of positive solutions. 
In contrast to these works, we impose invariance only under the Kelvin transform, rather than under a larger group of conformal transformations. 
Although this weaker symmetry does not restore compactness of the critical Sobolev embedding, 
it allows us to analyze possible concentration bubbles by combining a profile decomposition argument with suitable energy estimates. 
Moreover, unlike the constant-coefficient settings considered in these works, we allow the coupling coefficient $\beta(x)$ to depend on $x$, subject to the Kelvin-invariance condition $(\beta2)$.

We now introduce the Kelvin transform $K$ used in this paper.
For $u\in C_0^\infty(\R^4\setminus\{0\})$, define
\begin{align}
(Ku)(x) = \frac{1}{|x|^{2}} u \left( -\frac{x}{|x|^2} \right).
\end{align}
Since $\|Ku\|_{\mathcal{D}^{1,2}}=\|u\|_{\mathcal{D}^{1,2}}$ and 
$C_0^\infty(\R^4 \setminus\{0\})$ is dense in $\mathcal{D}^{1,2}(\R^4)$,
$K$ extends uniquely to a linear isometry on $\mathcal{D}^{1,2}(\R^4)$. 
Note that the Kelvin transform is typically defined as
\begin{equation*}
	(\widetilde{K}u)(x) = \frac{1}{|x|^{2}} u \left( \frac{x}{|x|^2} \right).
\end{equation*}
However, the inversion $x \mapsto \frac{x}{\vert{}x\vert{}^2}$ fixes every point on the unit sphere, whereas $x \mapsto -\frac{x}{\vert{}x\vert{}^2}$ has no fixed points, which motivates our definition of $K$ above. 
We define the subspace of functions invariant under the Kelvin transform by 
\begin{equation*}
\mathcal{D}_K := \{\bm{u} \in \mathcal{D} \mid Ku_1 = u_1, Ku_2=u_2\}.
\end{equation*}
A function $\bm{u} \in \mathcal{D}_K$ is called Kelvin-invariant. Furthermore, we say that $\bm{u}$ is positive if both of its components are positive.

We consider the following minimization problem:
\begin{align*}	
c_K&:= \inf_{\bm{u}\in\mathcal{N}_K} I(\bm{u}),\\
\mathcal{N}_K &:= \{\bm{u} \in \mathcal{D}_K \mid I^\prime(\bm{u})(u_1,0)=0, I^\prime(\bm{u})(0,u_2)=0, u_1\not\equiv 0, u_2\not\equiv 0\}.
\end{align*}
As shown in Lemma \ref{lem 2.3}, $\mathcal{N}_K$ is a natural constraint for $I$.
Thus, if $c_K$ is attained, then the corresponding minimizer is a critical point of $I$.
Our first main result is the following.

\begin{thm}
\label{thm 1.1}
Suppose that $\mu_1\leq \mu_2$ and ($\beta$1)--($\beta$2) hold. 
Set $t=\frac{\mu_1}{\mu_2}\in (0,1]$.
If
\begin{equation}
	\label{1.3}
	\int_{\R^4} \frac{\beta(x)}{(1+|x|^2)^4 } \, dx <  \frac{t}{6(1 + \sqrt{2+2t})} \cdot \pi^2\mu_2,
\end{equation}
then $c_K$ is attained and \eqref{P} has a Kelvin-invariant positive solution.
\end{thm}

Theorem \ref{thm 1.1} also covers the case where $\beta$ is constant. 
Note that $t  \mapsto \frac{t}{6(1 + \sqrt{2+2t})}$ is strictly increasing and satisfies 
\begin{equation*}
	\lim_{t\to 0^+}\frac{t}{6(1 + \sqrt{2+2t})}=0 \quad \text{and} \quad
	\lim_{t\to 1}\frac{t}{6(1 + \sqrt{2+2t})} = \frac{1}{18}.
\end{equation*}
We also obtain an existence result under a different weighted integral condition on $\beta$.

\begin{thm}
\label{thm 1.2}
Suppose that $\mu_1\leq \mu_2$ and ($\beta$1)--($\beta$2) hold. 
If
\begin{equation}
\label{1.4}
 \int_{\R^4} \frac{\beta(x)}{|x|^4} \, dx <  \pi^2 \mu_2,
\end{equation}
then $c_K$ is attained and \eqref{P} has a Kelvin-invariant positive solution.
\end{thm}

Theorem \ref{thm 1.2} does not cover the case where $\beta$ is constant.

\begin{rmk}
Neither condition \eqref{1.3} nor condition \eqref{1.4} implies the other. 
Indeed, for every $t\in(0,1]$, the relative sizes of the following two weights depend on $|x|$:
\begin{equation*}
	\frac{6(1+\sqrt{2+2t})}{t(1+|x|^2)^4} \quad \text{and} \quad \frac{1}{|x|^4}.
\end{equation*}
The former is larger near $|x|=1$, while the latter is larger near $|x|=0$ and for sufficiently large $|x|$.
Hence, which condition is more restrictive depends on $\beta$.
\end{rmk}

The proofs of Theorems \ref{thm 1.1} and \ref{thm 1.2} are based on a profile decomposition argument. 
Using the Kelvin invariance and energy estimates for nontrivial profiles, 
we show that every Palais--Smale sequence for $I$ contained in $\mathcal{N}_K$ at the level $c_K$ has a strongly convergent subsequence, provided that
\begin{equation*}
	c_K < \frac{1}{4}\left(\frac{1}{\mu_1} + \frac{2}{\mu_2}\right)S^2,
\end{equation*}
where $S$ is the best Sobolev constant given by
\begin{equation*}
	S=\inf_{u\in\mathcal{D}^{1,2}(\R^4)\setminus \{0\}}\frac{\|u\|^2_{\mathcal{D}^{1,2}}}{\|u\|^2_{4}}.
\end{equation*}
To obtain this estimate, we construct explicit test functions.
We define the Talenti function $U$ by
\begin{equation*}
 U(x) = \frac{2^{\frac{3}{2}}}{1+|x|^2}.
\end{equation*}
Also, we employ a family of Talenti functions $(U_{\varepsilon})_{\varepsilon>0}$ defined by 
\begin{equation}
\label{1.5}
	U_\varepsilon(x) = \varepsilon^{-1}U(\varepsilon^{-1}x)
	= \frac{2^{\frac{3}{2}}\varepsilon}{\varepsilon^2 + |x|^2}. 
\end{equation}
Note that $U_\varepsilon$ attains $S$ and is a positive solution of $-\Delta u = u^{3}$ in $\R^4$ (see \cite{Talenti}). Under the assumptions of Theorem \ref{thm 1.1}, we use $(U_1,U_1)$ as a test function, which is the simplest choice. 
By contrast, to prove Theorem \ref{thm 1.2}, we employ $(U_\varepsilon + U_{\frac{1}{\varepsilon}},U_1)$ instead. 
This type of test function is known to induce an interaction that lowers the energy of the functional (see, e.g. \cite{Adachi, Bahri--Li, Hirata08, Sato23, Sato--Suzuki}).
The conditions \eqref{1.3} and \eqref{1.4} are sharp with respect to the estimates obtained from the test functions $(U_1,U_1)$ and $(U_\varepsilon + U_{\frac{1}{\varepsilon}},U_1)$.

Furthermore, if $\beta(x)$ is a radial function, 
we obtain the following existence result for radially symmetric solutions by essentially the same argument as in the proofs of
Theorems \ref{thm 1.1} and \ref{thm 1.2}.
Indeed, the test functions $(U_1,U_1)$ and $(U_\varepsilon + U_{\frac{1}{\varepsilon}},U_1)$ are radially symmetric, 
so the same energy estimates also hold in the radially symmetric class.
We consider the following minimization problem:
\begin{align*}
	c_{K, {\rm rad}}&:= \inf_{\bm{u}\in\mathcal{N}_K \cap \mathcal{D}_{\rm rad}} I(\bm{u}),\\
	\mathcal{D}_{\rm rad} &:= \{\bm{u} \in \mathcal{D} \mid u_1, u_2 \text{ are radial functions}\}.
\end{align*}
If $c_{K,{\rm rad}}$ is attained, then any minimizer is also a critical point of $I$.

\begin{thm}
\label{thm 1.4}
Suppose that the assumptions of Theorem \ref{thm 1.1} or Theorem \ref{thm 1.2} hold. If $\beta(x)$ is a radial function, 
then $c_{K,{\rm rad}}$ is attained and \eqref{P} has a Kelvin-invariant and radial positive solution.
\end{thm}

\begin{rmk}
When $\beta$ is a positive constant, Chen, Medina, and Pistoia \cite[Theorem 1.1]{Chen--Medina--Pistoia} constructed nonradial solutions, 
while Clapp and Pistoia \cite[Theorem 1.1]{Clapp--Pistoia} obtained positive solutions under stronger symmetry assumptions.
On the other hand, Theorems \ref{thm 1.1}, \ref{thm 1.2} and \ref{thm 1.4} establish the existence of least-energy solutions in the Kelvin-invariant class 
and the radially symmetric Kelvin-invariant class, respectively, allowing $\beta$ to be nonconstant.
\end{rmk}

Moreover, we study the asymptotic behavior of the minimizers of $c_K$ and $c_{K,{\rm rad}}$ as $\|\beta\|_\infty \to 0$. 
Note that, for $\mu > 0$, $\mu^{-\frac{1}{2}} U$ is the unique positive Kelvin-invariant solution of the following equation.
\begin{equation}
\label{1.6}
-\Delta u = \mu u^3 \quad \text{in } \R^4. 
\end{equation}

\begin{thm}
\label{thm 1.6}
Suppose that ($\beta$1)--($\beta$2) hold. 
Let $\bm{w}_\beta$ be a positive minimizer of $c_K$ obtained in Theorem \ref{thm 1.1} or \ref{thm 1.2}, or a positive minimizer
of $c_{K,{\rm rad}}$ obtained in Theorem \ref{thm 1.4}. 
Then, as $\|\beta\|_\infty \to 0$, the following holds:
\begin{equation*}
	\bm{w}_\beta \to (\mu_1^{-\frac{1}{2}}U, \mu_2^{-\frac{1}{2}}U) \quad \text{ strongly in $\mathcal{D}_K$}.
\end{equation*}
\end{thm}
We emphasize that the convergence in Theorem \ref{thm 1.6} holds without passing to a subsequence.

This paper is organized as follows.
In Section 2, we provide some preliminaries.
In Section 3, we derive energy estimates.
In Section 4, we complete the proofs of Theorems \ref{thm 1.1}, \ref{thm 1.2}, and \ref{thm 1.4}.
Finally, in Section 5, we investigate the asymptotic behavior of minimizers of $c_K$ and $c_{K,{\rm rad}}$.

\section{Preliminaries}
In this section, we recall some properties of the Kelvin transform and prepare preliminary results for the subsequent analysis.

\begin{lem}
\label{lem 2.1}
The Kelvin transform $K$ on $\mathcal{D}^{1,2}(\R^4)$ satisfies the following properties:
\begin{itemize}
\item[(i)] $K^2 = id$.
\item[(ii)] $\|Ku\|_{\mathcal{D}^{1,2}(\R^4)} = \|u\|_{\mathcal{D}^{1,2}(\R^4)}$ for all $u \in \mathcal{D}^{1,2}(\R^4)$.
\item[(iii)] $\|Ku\|_{4} = \|u\|_{4}$ for all $u \in L^{4}(\R^4)$.
\item[(iv)] If ($\beta$2) holds, we have $\displaystyle B(Ku_1,Ku_2) = B(u_1,u_2)$ for all $\bm{u} \in \mathcal{D}$.
\end{itemize}
\end{lem}

\begin{proof}
Properties (i)--(iii) are standard; their proofs can be found in
\cite[Lemma 3.1]{Sato--Suzuki}.
It remains to prove (iv).
For $\bm{u} \in C_0^\infty(\R^N \setminus \{0\})^2$, we have 
\begin{equation*}
B(Ku_1,Ku_2) = \int_{\R^4} \beta(x)(Ku_1)^2 (Ku_2)^2 \, dx = \int_{\R^4} \beta(x) \frac{1}{|x|^8} u_1\left(- \frac{x}{|x|^2}\right)^2 u_2\left(- \frac{x}{|x|^2}\right)^2 \, dx.
\end{equation*}
By the change of variables $x = - \frac{y}{|y|^2}$, whose absolute Jacobian determinant is $|y|^{-8}$, and ($\beta$2), we have 
\begin{equation*}
 B(Ku_1,Ku_2) = \int_{\R^4} \beta\left(- \frac{y}{|y|^2}\right) u_1^2 u_2^2 \, dy = \int_{\R^4} \beta(y) u_1^2 u_2^2 \, dy.
\end{equation*}
Since $C_0^\infty(\R^N \setminus \{0\})^2$ is dense in $\mathcal{D}$, (iv) holds. 
\end{proof}

By Lemma \ref{lem 2.1}, if ($\beta$2) holds, then we have
\begin{equation*}
I(Ku_1,Ku_2) = I(u_1,u_2) \quad \text{for all } \bm{u} \in \mathcal{D},
\end{equation*}
that is, $I$ is invariant under the Kelvin transform.

Next, we establish some basic variational properties of $I$ on the Nehari set $\mathcal{N}_K$.
In particular, we derive an alternative characterization of $c_K$, which will be useful for the energy estimates in Section 3.
For any $\bm{u} \in \mathcal{N}_K$, we have
\begin{equation}
\label{2.1}
\|u_1\|^2_{\mathcal{D}^{1,2}} = \mu_1\|u_1\|_{4}^4 - B(\bm{u}), \quad \|u_2\|^2_{\mathcal{D}^{1,2}} = \mu_2\|u_2\|_{4}^4 - B(\bm{u}).
\end{equation}
Consequently, $I(\bm{u})$ can be written as
\begin{equation}\label{2.2}
	I(\bm{u}) = \frac{1}{4} \left(\|u_1\|^2_{\mathcal{D}^{1,2}} + \|u_2\|^2_{\mathcal{D}^{1,2}}\right) > 0 
	\quad \text{for all } \bm{u} \in \mathcal{N}_K.
\end{equation}
Moreover, \eqref{2.1} implies $\mu_i\|u_i\|_{4}^4> B(\bm{u})$ for $i=1,2$. 
Thus we have
\begin{equation}\label{2.0}
	\mu_1\mu_2\|u_1\|_{4}^4\|u_2\|_{4}^4 - B(\bm{u})^2> 0 
	\quad \text{for all } \bm{u} \in \mathcal{N}_K.
\end{equation}
We begin with the following uniform lower bound.

\begin{lem}
\label{lem 2.2}
Suppose that ($\beta$1) holds. Then, there exists a constant $\delta > 0$ such that for all $\bm{u} \in \mathcal{N}_K$, $\|u_i\|_{4} \ge \delta$ for $i = 1, 2$.
\end{lem}

\begin{proof}
Let $(u_1,u_2) \in \mathcal{N}_K$. By \eqref{2.1} and the Sobolev inequality, we have 
\begin{equation*}
\|u_i\|_{4}^2 \le S^{-1} \|u_i\|_{\mathcal{D}^{1,2}}^2 \le S^{-1} \mu_i \|u_i\|_{4}^4.
\end{equation*}
Therefore, $\|u_i\|_{4} \ge \delta$ holds with $\delta = S^{\frac{1}{2}} \min\{ \mu_1^{-\frac{1}{2}}, \mu_2^{-\frac{1}{2}}\}$.
\end{proof}

The following lemma shows that $\mathcal{N}_K$ is a natural constraint.

\begin{lem}
\label{lem 2.3}
Suppose that ($\beta$1) and ($\beta$2) hold. 
Then, $\mathcal{N}_K$ is a natural constraint of $I$, that is, critical points of $I|_{\mathcal{N}_K}$ are critical points of $I$.
\end{lem}

\begin{proof}
First, we show that $\mathcal{N}_K$ is a $C^1$-manifold of codimension 2. 
Define $G_1,G_2:\mathcal{D}_K\to\R$ by
\begin{equation*}
	G_1(\bm{u}) = I^\prime(\bm{u})(u_1,0), \quad G_2(\bm{u}) = I^\prime(\bm{u})(0,u_2).
\end{equation*}
Then, $\mathcal{N}_K$ can be written as
\begin{equation*}
	\mathcal{N}_K = \{\bm{u} \in \mathcal{D}_K\mid G_1(\bm{u})=0, G_2(\bm{u})=0, u_1\not= 0, u_2\not= 0 \}.
\end{equation*}
Let $\bm{u}\in\mathcal{N}_K$.
We show that $G_1^\prime (\bm{u}), G_2^\prime (\bm{u})$ are linearly independent. 
Suppose that the following identity holds for some $c_1, c_2 \in \R$:
\begin{equation*}
c_1 G^\prime_1 (\bm{u}) \bm{\varphi} + c_2 G^\prime_2 (\bm{u})\bm{\varphi} = 0 \quad \text{for all } \bm{\varphi} \in \mathcal{D}_K.
\end{equation*}
Setting $\bm{\varphi}=(u_1,0)$ and $(0,u_2)$, we obtain
\begin{align*}
\begin{pmatrix}
-2\mu_1\|u_1\|_4^4 & 2B(\bm{u})\\
2B(\bm{u}) &-2\mu_2\|u_2\|_4^4
\end{pmatrix}
\begin{pmatrix}
c_1\\
c_2
\end{pmatrix}
=
\begin{pmatrix}
0\\
0
\end{pmatrix}.
\end{align*}
By \eqref{2.0}, the determinant of the above matrix is positive, which implies $(c_1,c_2) = (0,0)$. 
Hence, $G_1^\prime (\bm{u}), G_2^\prime (\bm{u})$ are linearly independent. 
Therefore, $\mathcal{N}_K$ is a $C^1$-manifold of codimension 2. Let $\bm{u}$ be a critical point of $I|_{\mathcal{N}_K}$. 
Then, by the method of Lagrange multipliers (see \cite[Theorem 6.2.1]{Costa}), there exist $\lambda_1, \lambda_2 \in \R$ such that 
\begin{equation*}
I^\prime(\bm{u})\bm{\varphi} = \lambda_1 G_1^\prime (\bm{u}) \bm{\varphi} + \lambda_2 G_2^\prime (\bm{u}) \bm{\varphi} 
\quad \text{for all } \bm{\varphi} \in \mathcal{D}_K.
\end{equation*}
Setting $\bm{\varphi}=(u_1,0)$ and $(0,u_2)$ and using $G_1(\bm{u})=G_2(\bm{u})=0$, we obtain
\begin{align*}
\begin{pmatrix}
-2\mu_1\|u_1\|_4^4 & 2B(\bm{u})\\
2B(\bm{u}) &-2\mu_2\|u_2\|_4^4
\end{pmatrix}
\begin{pmatrix}
\lambda_1\\
\lambda_2
\end{pmatrix}
=
\begin{pmatrix}
0\\
0
\end{pmatrix}.
\end{align*}
Since the determinant is positive, we have $\lambda_1=\lambda_2=0$. 
Thus, $\bm{u}$ is a critical point of $I|_{\mathcal{D}_K}$. Moreover, since $I$ is invariant under the Kelvin transform, Palais' principle of symmetric criticality \cite{Palais} yields that $\bm{u}$ is a critical point of $I$. Therefore, $\mathcal{N}_K$ is a natural constraint of $I$.
\end{proof}

Set
\begin{equation*}
\mathcal{X}_K := \{\bm{u}\in\mathcal{D}_K\mid \mu_1 \mu_2\|u_1\|_4^4 \|u_2\|_4^4 - B(\bm{u})^2 >0\}.
\end{equation*}
Then, \eqref{2.0} implies $\mathcal{N}_K \subset \mathcal{X}_K$. 
We now establish the fibering map property of the functional $I$ on $\mathcal{X}_K$. 

\begin{lem}
\label{lem 2.4}
Suppose that ($\beta$1) holds.
For $\bm{u} \in \mathcal{X}_K$, there exists a unique pair $(s(\bm{u}),t(\bm{u})) \in (0,\infty)^2$ such that $(s(\bm{u})u_1,t(\bm{u})u_2) \in \mathcal{N}_K$ and 
\begin{align*}
\begin{pmatrix}
s(\bm{u})^2\\
t(\bm{u})^2
\end{pmatrix}
&=\frac{1}{\mu_1 \mu_2 \|u_1\|_4^4 \|u_2\|_4^4 - B(\bm{u})^2}
\begin{pmatrix}
  \mu_2\|u_2\|_4^4 &B(\bm{u})\\
  B(\bm{u}) &\mu_1\|u_1\|_4^4
\end{pmatrix}
\begin{pmatrix}
\|u_1\|^2_{\mathcal{D}^{1,2}}\\
\|u_2\|^2_{\mathcal{D}^{1,2}}
\end{pmatrix}.
\end{align*}
\end{lem}

\begin{proof} 
Let $\bm{u} \in \mathcal{X}_K$. We define $f : (0,\infty)^2 \to \mathbb{R}$ by
\begin{align*}
f(s,t) &= I(su_1,tu_2) = \frac{1}{2}s^2 \|u_1\|^2_{\mathcal{D}^{1,2}} + \frac{1}{2}t^2 \|u_2\|^2_{\mathcal{D}^{1,2}} - \frac{\mu_1}{4}s^4\|u_1\|_4^4 - \frac{\mu_2}{4}t^4\|u_2\|_4^4 + \frac{1}{2}s^2 t^2 B(\bm{u}).
\end{align*}
Then 
\begin{equation*}
	\nabla f(s,t) = \begin{pmatrix}
	I'(su_1,tu_2)(u_1,0) \\ I'(su_1,tu_2)(0,u_2)
	\end{pmatrix}
	=\begin{pmatrix}
	s\|u_1\|^2_{\mathcal{D}^{1,2}} - s^3\|u_1\|_4^4 + st^2B(\bm{u}) \\ 
	t\|u_2\|^2_{\mathcal{D}^{1,2}} - t^3\|u_2\|_4^4 + s^2tB(\bm{u})
	\end{pmatrix}.
\end{equation*}
We see that $\nabla f(s,t) = \bm{0}$ is equivalent to the following equation:
\begin{align*}
\begin{pmatrix}
  \mu_1\|u_1\|_4^4 &-B(\bm{u})\\
  -B(\bm{u}) &\mu_2\|u_2\|_4^4
\end{pmatrix}
\begin{pmatrix}
s^2\\
t^2
\end{pmatrix}
=
\begin{pmatrix}
\|u_1\|^2_{\mathcal{D}^{1,2}}\\
\|u_2\|^2_{\mathcal{D}^{1,2}}
\end{pmatrix}.
\end{align*}
Since the determinant of the above matrix is positive, $\nabla f(s,t) = \bm{0}$ has the unique solution $(s(\bm{u}),t(\bm{u}))$ given by 
\begin{align*}
\begin{pmatrix}
s(\bm{u})^2\\
t(\bm{u})^2
\end{pmatrix}
&=\frac{1}{\mu_1 \mu_2 \|u_1\|_4^4 \|u_2\|_4^4-B(\bm{u})^2}
\begin{pmatrix}
  \mu_2\|u_2\|_4^4 &B(\bm{u})\\
  B(\bm{u}) &\mu_1\|u_1\|_4^4
\end{pmatrix}
\begin{pmatrix}
\|u_1\|^2_{\mathcal{D}^{1,2}}\\
\|u_2\|^2_{\mathcal{D}^{1,2}}
\end{pmatrix}.
\end{align*}
Hence, $(s(\bm{u})u_1,t(\bm{u})u_2) \in \mathcal{N}_K$. 
\end{proof}

We next construct a (PS) sequence at the level $c_K$.
Recall that, for $(\bm{u}_k) \subset \mathcal{D}$ and $d \in \R$, we say that $(\bm{u}_k)$ is a (PS)$_d$ sequence for $I$ if $I(\bm{u}_k) \to d$ and $I'(\bm{u}_k) \to 0$ as $k \to \infty$.

\begin{lem}
\label{lem 2.3.5}
Suppose that ($\beta$1)--($\beta$2) hold.
Then there exists a (PS)$_{c_K}$ sequence $(\bm{u}_k)\subset\mathcal{N}_K$ for $I$.
\end{lem}

\begin{proof}
Let $(\bm{v}_k)\subset\mathcal{N}_K$ be a minimizing sequence for $c_K$.
By Lemma \ref{lem 2.2}, $\mathcal{N}_K$ is a closed subset of $\mathcal{D}_K$.
Thus, by Ekeland's variational principle (cf. \cite[Theorem 1.1]{Ghoussoub}),
there exists $(\bm{u}_k)\subset\mathcal{N}_K$ such that $I(\bm{u}_k)\to c_K$ and
\begin{equation}\label{eq:ekeland} 
	I(\bm{v})\geq I(\bm{u}_k) -\frac{1}{k}\|\bm{v}-\bm{u}_k\|_{\mathcal{D}}
	\quad\text{for all }\bm{v}\in\mathcal{N}_K.
\end{equation}
By \eqref{2.2}, $(\bm{u}_k)$ is bounded in $\mathcal{D}$.
For $\bm{u}\in\mathcal{X}_K$, set
$
	\Phi(\bm{u}) = (s(\bm{u})u_1,t(\bm{u})u_2),
$
where $s(\bm{u})$ and $t(\bm{u})$ are given in Lemma \ref{lem 2.4}.
For any $\bm{h}\in\mathcal{D}_K$ and sufficiently small $\tau>0$, we have
$\Phi(\bm{u}_k+\tau\bm{h})\in\mathcal{N}_K$.
Hence, by \eqref{eq:ekeland},
\begin{equation*}
	I(\Phi(\bm{u}_k+\tau\bm{h}))-I(\bm{u}_k) 
	\geq -\frac{1}{k} \|\Phi(\bm{u}_k+\tau\bm{h})-\bm{u}_k\|_{\mathcal{D}}.
\end{equation*}
Dividing by $\tau$ and letting $\tau\to0^+$, we obtain
\begin{equation*}
	I'(\bm{u}_k)\Phi'(\bm{u}_k)\bm{h}
	\geq -\frac{1}{k}\|\Phi'(\bm{u}_k)\bm{h}\|_{\mathcal{D}}.
\end{equation*}
Replacing $\bm{h}$ by $-\bm{h}$ yields
\begin{equation*}
	|I'(\bm{u}_k)\Phi'(\bm{u}_k)\bm{h}|
	\leq \frac{1}{k}\|\Phi'(\bm{u}_k)\bm{h}\|_{\mathcal{D}}.
\end{equation*}
Since $\bm{u}_k\in\mathcal{N}_K$, we have $s(\bm{u}_k)=t(\bm{u}_k)=1$.
Hence
\begin{equation*}
	\Phi'(\bm{u}_k)\bm{h}=([s'(\bm{u}_k)h]u_{1,k}, [t'(\bm{u}_k)h]u_{2,k}) + \bm{h}.
\end{equation*}
Using $I'(\bm{u}_k)(u_{1,k},0) = I'(\bm{u}_k)(0,u_{2,k}) = 0$, we obtain
\begin{equation*}
	I'(\bm{u}_k)\Phi'(\bm{u}_k)\bm{h} = I'(\bm{u}_k)\bm{h}.
\end{equation*}
Moreover, by \eqref{2.1}, the boundedness of $(\bm{u}_k)$, and Lemmas \ref{lem 2.2} and \ref{lem 2.4}, 
the derivatives $\Phi'(\bm{u}_k)$ are uniformly bounded.
Consequently,
\begin{equation*}
	|I'(\bm{u}_k)\bm{h}| \leq \frac{C}{k}\|\bm{h}\|_{\mathcal{D}}
	\quad\text{for all }\bm{h}\in\mathcal{D}_K.
\end{equation*}
Thus, $I'(\bm{u}_k)\to0$ in $\mathcal{D}_K^*$.
Since $I$ is invariant under the Kelvin transform, we also have $I'(\bm{u}_k) \to 0$ in $\mathcal{D}^*$.
Therefore, $(\bm{u}_k)$ is a (PS)$_{c_K}$ sequence for $I$.
\end{proof}

Finally, we provide an alternative characterization of $c_K$.
By Lemma \ref{lem 2.4}, for $\bm{u}\in \mathcal{X}_K$, we have $(s(\bm{u})u_1,t(\bm{u})u_2) \in \mathcal{N}_K$.
Hence by \eqref{2.2} and Lemma \ref{lem 2.4}, we obtain
\begin{align*}
	I(s(\bm{u})u_1, t(\bm{u})u_2) & = \frac{1}{4}(s(\bm{u})^2\|u_1\|^2_{\mathcal{D}^{1,2}} + t(\bm{u})^2\|u_2\|^2_{\mathcal{D}^{1,2}}) \\
	& = \frac{1}{4}\frac{\mu_2\|u_2\|^4_4 \|u_1\|^4_{\mathcal{D}^{1,2}} + 2B(\bm{u})\|u_1\|_{\mathcal{D}^{1,2}}^2 \|u_2\|_{\mathcal{D}^{1,2}}^2 + \mu_1\|u_1\|^4_4 \|u_2\|^4_{\mathcal{D}^{1,2}}}{\mu_1\mu_2\|u_1\|_4^4\|u_2\|_4^4-B(\bm{u})^2}.
\end{align*}
Define $J : \mathcal{X}_K \to \R$ by
\begin{align*}
	J(\bm{u}) = \frac{1}{4}\frac{\mu_2\|u_2\|^4_4 \|u_1\|^4_{\mathcal{D}^{1,2}} + 2B(\bm{u})\|u_1\|_{\mathcal{D}^{1,2}}^2 \|u_2\|_{\mathcal{D}^{1,2}}^2 + \mu_1\|u_1\|^4_4 \|u_2\|^4_{\mathcal{D}^{1,2}}}{\mu_1\mu_2\|u_1\|_4^4\|u_2\|_4^4-B(\bm{u})^2}.
\end{align*}
The following proposition provides another characterization of $c_K$.

\begin{prop}\label{prop 2.5}
Suppose that ($\beta$1) and ($\beta$2) hold. Then,
\begin{equation}
\label{2.3}
	c_K = \inf_{\bm{u}\in\mathcal{X}_K}J(\bm{u}).
\end{equation}
\end{prop}

\begin{proof}
Let $\bm{u} \in \mathcal{X}_K$. 
Then, we have 
\begin{equation*}
c_K \le I(s(\bm{u})u_1, t(\bm{u})u_2)=J(\bm{u}).
\end{equation*}
On the other hand, let $\bm{u} \in \mathcal{N}_K$. 
Then the uniqueness in Lemma \ref{lem 2.4} yields $s(\bm{u})=t(\bm{u})=1$. 
Since $\mathcal{N}_K\subset\mathcal{X}_K$, we obtain
\begin{equation*}
c_K = \inf_{\bm{u}\in\mathcal{N}_K}I(s(\bm{u})u_1,t(\bm{u})u_2) \ge \inf_{\bm{u}\in\mathcal{X}_K}I(s(\bm{u})u_1,t(\bm{u})u_2)\ = \inf_{\bm{u}\in\mathcal{X}_K}J(\bm{u}).
\end{equation*}
Therefore, the proof is complete.
\end{proof}

\section{Energy estimate}
In this section, under the assumptions of Theorem \ref{thm 1.1} or Theorem \ref{thm 1.2}, 
we prove the energy estimate
\begin{equation}
\label{3.1}
	c_K < \frac{1}{4}\left(\frac{1}{\mu_1} + \frac{2}{\mu_2}\right)S^2.
\end{equation}
Recall that $U_\varepsilon$, defined by \eqref{1.5}, is a positive solution of 
$-\Delta u = u^{3}$ in $\R^4$ and attains the best Sobolev constant $S$.
Thus $U_{\varepsilon}$ satisfies
\begin{align}
\label{3.2}
\int_{\R^4} \nabla U_{\varepsilon} \cdot \nabla \varphi \, dx 
& = \int_{\R^4} U_{\varepsilon}^{3}\varphi \, dx \quad \text{for all } \varphi \in \mathcal{D}^{1,2}(\R^4), \\
\label{3.3}
\|U_\varepsilon\|_{\mathcal{D}^{1,2}}^2 & = \|U_\varepsilon\|^4_4 =S^2.
\end{align}
Moreover, a direct calculation shows that 
\begin{equation}\label{3.4}
	KU_\varepsilon = U_{\frac{1}{\varepsilon}} \quad \text{ for } \varepsilon>0.
\end{equation}

We first prove \eqref{3.1} under the assumptions of Theorem \ref{thm 1.1}.

\begin{prop}
\label{prop 3.1}
Suppose that $\mu_1\leq \mu_2$ and ($\beta$1)--($\beta$2) hold. 
Then \eqref{1.3} implies \eqref{3.1}. 
\end{prop}

\begin{proof}
To use $(U_1,U_1)$ as a test function, we first confirm that $(U_1,U_1) \in \mathcal{X}_K$.
By \eqref{3.4} with $\varepsilon = 1$, we see that $(U_1,U_1) \in \mathcal{D}_K$. 
Since $U_1(x)=\frac{2^{\frac{3}{2}}}{1+|x|^2}$, we have
\begin{equation*}
	B(U_1,U_1)=\int_{\R^4} \beta(x) U_1^4 \, dx 
	= 64 \int_{\R^4} \frac{\beta(x)}{(1+|x|^2)^4} \, dx.
\end{equation*}
From \eqref{1.3}, with $t=\frac{\mu_1}{\mu_2}$, we obtain
\begin{equation*}
	B(U_1,U_1)< \frac{32}{3} \frac{t}{1+\sqrt{2+2t}} \cdot \pi^2 \mu_2 
	< \frac{32}{3} \sqrt{t} \pi^2 \mu_2 
	= \frac{32}{3} \sqrt{\mu_1\mu_2} \pi^2.  
\end{equation*}
Since $S^2 = \|U_1\|_4^4 = \frac{32}{3} \pi^2$ (see Remark \ref{rmk 3.4}), we have
\begin{equation*}
	\mu_1 \mu_2 \|U_1\|_4^8  - B(U_1,U_1)^2  > 0.
\end{equation*}
Hence, $(U_1,U_1) \in \mathcal{X}_K$.
By Proposition \ref{prop 2.5}, we have
\begin{equation}
\label{3.5}
\begin{aligned}
c_K &\le J(U_1,U_1)\\
&= \frac{1}{4}\frac{\mu_2\|U_1\|^4_4 \|U_1\|^4_{\mathcal{D}^{1,2}} + 2B(U_1,U_1)\|U_1\|_{\mathcal{D}^{1,2}}^2 \|U_1\|_{\mathcal{D}^{1,2}}^2 + \mu_1\|U_1\|^4_4 \|U_1\|^4_{\mathcal{D}^{1,2}}}{\mu_1\mu_2\|U_1\|_4^4\|U_1\|_4^4-B(U_1,U_1)^2}\\
&= \frac{1}{4}\frac{\mu_2 S^2 \cdot S^4 + 2B(U_1,U_1) \cdot S^2 \cdot S^2 + \mu_1 S^2 \cdot S^4}{\mu_1 \mu_2 S^2 \cdot S^2 - B(U_1,U_1)^2}\\
&= \frac{1}{4} \frac{(\mu_1 + \mu_2)S^6 + 2B(U_1,U_1)S^4}{\mu_1 \mu_2 S^4 - B(U_1,U_1)^2}.
\end{aligned}
\end{equation}
To complete the proof, we observe that
\begin{align*}
& \frac{1}{4} \frac{(\mu_1 + \mu_2)S^6 + 2B(U_1,U_1)S^4}{\mu_1 \mu_2 S^4 - B(U_1,U_1)^2} 
< \frac{1}{4}\left(\frac{1}{\mu_1} + \frac{2}{\mu_2}\right)S^2 \\
& \iff \left(\frac{1}{\mu_1} + \frac{2}{\mu_2}\right)B(U_1,U_1)^2 + 2S^2 B(U_1,U_1) - \mu_1S^4 < 0 \\
& \iff B(U_1,U_1) < \frac{t}{1+\sqrt{2+2t}} \cdot \mu_2 S^2.
\end{align*}
The last inequality follows from (1.3).
Combining this with \eqref{3.5}, we obtain \eqref{3.1}.
\end{proof}

Next, we establish \eqref{3.1} under the assumptions of Theorem 1.2.

\begin{prop}
\label{prop 3.2}
Suppose that $\mu_1\leq \mu_2$ and ($\beta$1)--($\beta$2) hold. 
Then \eqref{1.4} implies \eqref{3.1}.
\end{prop}

To prove Proposition 3.2, we first establish some asymptotic estimates as $\varepsilon \to 0^+$.

\begin{lem}
\label{lem 3.3}
As $\varepsilon \to 0^+$, the following hold:
\begin{itemize}
\item[(i)] $\displaystyle \int_{\R^4} U_\varepsilon^3 U_{\frac{1}{\varepsilon}} \, dx = A \varepsilon^2 + o(\varepsilon^2)$, where $A = 2^{\frac{3}{2}} \|U_1\|_3^3$.
\item[(ii)] $\displaystyle \int_{\R^4} U_\varepsilon^2 U_{\frac{1}{\varepsilon}}^2 \, dx = o(\varepsilon^2)$.
\end{itemize}
\end{lem}

\begin{proof}
First, we prove (i). By the change of variables $x=\varepsilon y$, we have
\begin{equation*}
\int_{\R^4}  U_\varepsilon^3 U_{\frac{1}{\varepsilon}} \, dx 
= \int_{\R^4} \varepsilon^3 U_1(\varepsilon^{-1}x)^3 \varepsilon U_1(\varepsilon x) \, dx = \varepsilon^2 \int_{\R^4} U_1(y)^3 U_1(\varepsilon^2 y) \, dy.
\end{equation*}
Since $U_1(\varepsilon^2 y) \to U_1(0) = 2^{\frac{3}{2}}$ pointwise as $\varepsilon \to 0$, we deduce
\begin{equation*}
\int_{\R^4} U_\varepsilon^3 U_{\frac{1}{\varepsilon}} \, dx = 2^{\frac{3}{2}} \left(\int_{\R^4} U_1^3 \, dx\right)\varepsilon^2 + o(\varepsilon^2). 
\end{equation*}
Hence, (i) holds.
Next, we prove (ii). Let $p, q > 1$ be H\"older conjugate exponents. A direct computation shows that 
\begin{equation*}
\|U_\varepsilon\|_r = \varepsilon^{\frac{4}{r} - 1} \|U_1\|_r \quad \text{for all } \varepsilon > 0 \text{ and } r > 1.
\end{equation*}
Combining this with the Hölder inequality, we have
\begin{align*}
\int_{\R^4} U_\varepsilon^2 U_{\frac{1}{\varepsilon}}^2 \, dx
&\le \|U_\varepsilon\|^2_{2p} \|U_{\frac{1}{\varepsilon}}\|^2_{2q}\\
&= \varepsilon^{2\left(\frac{2}{p} - 1\right)} \|U_1\|^2_{2p} \varepsilon^{-2\left(\frac{2}{q} - 1\right)} \|U_1\|^2_{2q} \\
&= \varepsilon^{4\left(\frac{1}{p} - \frac{1}{q}\right)} \|U_1\|^2_{2p} \|U_1\|^2_{2q}.
\end{align*}
Choosing $p = \frac{8}{7}$ and $q = 8$, we have $4\left(\frac{1}{p} - \frac{1}{q}\right) = 3$. Hence, we have
\begin{equation*}
\int_{\R^4} U_\varepsilon^2 U_{\frac{1}{\varepsilon}}^2 \, dx = O(\varepsilon^3) = o(\varepsilon^2).
\end{equation*}
Thus, (ii) holds. 
\end{proof}

\begin{rmk}
\label{rmk 3.4}
The best Sobolev constant $S$ and the constant $A$ in Lemma \ref{lem 3.3} can be computed explicitly as follows.
Indeed, for $p > 2$, integration in polar coordinates gives
\begin{equation*}
\int_{\R^4} \frac{1}{(1 + |x|^2)^p} \, dx = |\mathbb{S}^3| \int_{0}^\infty \frac{r^3}{(1 + r^2)^p} \, dr = \frac{\pi^2}{(p-1)(p-2)},
\end{equation*}
where $|\mathbb{S}^3|$ denotes the surface area of the unit sphere $\mathbb{S}^3$ in $\R^4$. 
Thus, we deduce
\begin{align*}
& S^2 = \|U_1\|_4^4 = \int_{\R^4} \frac{2^6}{(1 + |x|^2)^4} \, dx = \frac{32}{3} \pi^2, \\
& A =2^{\frac{3}{2}}\|U_1\|_3^3 = 2^{\frac{3}{2}} \int_{\R^4} \frac{2^{\frac{9}{2}}}{(1 + |x|^2)^3} \, dx = 32 \pi^2.
\end{align*}
\end{rmk}

\begin{lem}
\label{lem 3.5}
Suppose that $\frac{\beta(x)}{|x|^4} \in L^1(\R^4)$. Let $W_\varepsilon = U_\varepsilon + U_{\frac{1}{\varepsilon}}$. Then, as $\varepsilon \to 0^+$, the following hold:
\begin{itemize}
\item[(i)] $\|W_\varepsilon\|_{\mathcal{D}^{1,2}}^2 = 2S^2 + 2A\varepsilon^2 + o(\varepsilon^2)$.
\item[(ii)] $\|W_\varepsilon\|_4^4 = 2S^2 + 8A\varepsilon^2 + o(\varepsilon^2)$.
\item[(iii)] $\displaystyle B(W_\varepsilon,U_1) = 64 D \varepsilon^2 + o(\varepsilon^2)$.
\end{itemize}
Here $A = 2^{\frac{3}{2}}\|U_1\|_3^3$, $\displaystyle D = \int_{\R^4} \frac{\beta(x)}{|x|^4} \, dx$.
\end{lem}

\begin{proof}
Using \eqref{3.2} and \eqref{3.3}, we have
\begin{align*}
\|W_\varepsilon\|_{\mathcal{D}^{1,2}}^2 &= \int_{\R^4}|\nabla U_\varepsilon|^2 \, dx + \int_{\R^4}|\nabla U_{\frac{1}{\varepsilon}}|^2 \, dx + 2\int_{\R^4} \nabla U_\varepsilon \cdot \nabla U_{\frac{1}{\varepsilon}} \, dx\\
&= 2S^2 + 2\int_{\R^4} U_\varepsilon^3 U_{\frac{1}{\varepsilon}} \, dx.
\end{align*}
As $\varepsilon \to 0^+$, by Lemma \ref{lem 3.3} (i), we obtain
\begin{equation*}
\|W_\varepsilon\|_{\mathcal{D}^{1,2}}^2 = 2S^2 + 2A\varepsilon^2 + o(\varepsilon^2).
\end{equation*}
Hence, (i) holds. By the change of variables as in the proof of Lemma \ref{lem 2.1} and \eqref{3.4}, we have
\begin{equation*}
\int_{\R^4} U_\varepsilon^3 U_{\frac{1}{\varepsilon}} \, dx = \int_{\R^4} (KU_\varepsilon)^3 (KU_{\frac{1}{\varepsilon}}) \, dy  = \int_{\R^4} U_{\frac{1}{\varepsilon}}^3 U_\varepsilon \, dy.
\end{equation*}
Combining this with \eqref{3.3}, we obtain
\begin{align*}
\|W_\varepsilon\|_4^4 &= \int_{\R^4} U_\varepsilon^4 + 4 U_\varepsilon^3 U_{\frac{1}{\varepsilon}} + 6 U_\varepsilon^2 U_{\frac{1}{\varepsilon}}^2 + 4 U_\varepsilon U_{\frac{1}{\varepsilon}}^3 + U_{\frac{1}{\varepsilon}}^4 \,dx\\
&=2S^2 + 8 \int_{\R^4}U_\varepsilon^3 U_{\frac{1}{\varepsilon}} \, dx + 6 \int_{\R^4} U_\varepsilon^2 U_{\frac{1}{\varepsilon}}^2 \, dx.
\end{align*}
As $\varepsilon \to 0^+$, by Lemma \ref{lem 3.3} (i)--(ii), we have
\begin{equation*}
\|W_\varepsilon\|_4^4 = 2S^2 + 8A\varepsilon^2 + o(\varepsilon^2).
\end{equation*}
Thus, (ii) holds. 
Using the explicit form of the Talenti functions, we have
\begin{align*}
B(W_\varepsilon,U_1) &= \int_{\R^4} \beta(x) (U_\varepsilon + U_{\frac{1}{\varepsilon}})^2 U_1^2 \,dx\\
&= \int_{\R^4} \beta(x) \left(\frac{2^{\frac{3}{2}}\varepsilon}{\varepsilon^2 + |x|^2} + \frac{2^{\frac{3}{2}}\varepsilon}{1 + \varepsilon^2 |x|^2}\right)^2 \frac{8}{(1 + |x|^2)^2} \, dx\\
&= 64 \varepsilon^2 \int_{\R^4} \beta(x) \left(\frac{1}{\varepsilon^2 + |x|^2} + \frac{1}{1 + \varepsilon^2 |x|^2}\right)^2 \frac{1}{(1 + |x|^2)^2} \, dx.
\end{align*}
Futhermore, as $\varepsilon \to 0^+$, 
\begin{align*}
\left(\frac{1}{\varepsilon^2 + |x|^2} + \frac{1}{1 + \varepsilon^2 |x|^2}\right)^2 \frac{1}{(1 + |x|^2)^2} \to \left(\frac{1}{|x|^2} + 1 \right)^2 \frac{1}{(1 + |x|^2)^2} = \frac{1}{|x|^4} \quad \text{a.e. } x \in \R^4.
\end{align*}
Since $\frac{\beta(x)}{|x|^4} \in L^1(\R^4)$, by the dominated convergence theorem, we have
\begin{equation*}
B(W_\varepsilon,U_1) = 64 \varepsilon^2 \int_{\R^4} \frac{\beta(x)}{|x|^4} \, dx + o(\varepsilon^2).
\end{equation*}
Therefore, (iii) holds.
\end{proof}

\begin{proof}[Proof of Proposition \ref{prop 3.2}]
Set $W_\varepsilon = U_\varepsilon + U_{\frac{1}{\varepsilon}}$. 
To use $(W_\varepsilon,U_1)$ as a test function, 
we first confirm that $(W_\varepsilon,U_1) \in \mathcal{X}_K$ for sufficiently small $\varepsilon > 0$. 
By \eqref{3.4}, we have $(W_\varepsilon, U_1) \in \mathcal{D}_K$. 
As $\varepsilon \to 0^+$, by \eqref{3.3} and Lemma \ref{lem 3.5} (ii)--(iii), we obtain
\begin{align*}
	\mu_1 \mu_2 \|W_\varepsilon\|_4^4 \|U_1\|_4^4 - B(W_\varepsilon, U_1)^2 
	& = \mu_1 \mu_2 (2S^2 + 8A\varepsilon^2 + o(\varepsilon^2)) S^2 - (64D\varepsilon^2 + o(\varepsilon^2))^2\\
	& = 2\mu_1 \mu_2 S^4 + 8A\mu_1 \mu_2 S^2 \varepsilon^2 + o(\varepsilon^2) > 0.
\end{align*}
Thus, $(W_\varepsilon,U_1) \in \mathcal{X}_K$ for sufficiently small $\varepsilon > 0$. Furthermore, by using the expansion $(1+x)^{-1} = 1 - x +o(|x|)$ as $|x| \to 0$, it follows that
\begin{equation}
\label{3.7}
\begin{aligned}
\left(\frac{\mu_1 \mu_2 \|W_\varepsilon\|_4^4 \|U_1\|_4^4 - B(W_\varepsilon, U_1)^2}{2\mu_1 \mu_2 S^4}\right)^{-1}&= \left(1 + \frac{4A}{S^2} \varepsilon^2 + o(\varepsilon^2) \right)^{-1}\\
&= 1 - \frac{4A}{S^2} \varepsilon^2 + o(\varepsilon^2).
\end{aligned}
\end{equation}
As $\varepsilon \to 0^+$, by \eqref{3.3} and Lemma \ref{lem 3.5} (i)--(iii), we deduce
\begin{align*}
	& \mu_2\|U_1\|^4_4 \|W_\varepsilon\|^4_{\mathcal{D}^{1,2}} 
	  + 2B(W_\varepsilon,U_1)\|W_\varepsilon\|_{\mathcal{D}^{1,2}}^2 \|U_1\|_{\mathcal{D}^{1,2}}^2 
	  + \mu_1\|W_\varepsilon\|^4_4 \|U_1\|^4_{\mathcal{D}^{1,2}} \\
	& = \mu_2 S^2 (2S^2 + 2A\varepsilon^2 + o(\varepsilon^2))^2 + 2 (64D\varepsilon^2 + o(\varepsilon^2)) (2S^2 + 2A\varepsilon^2 
	  + o(\varepsilon^2)) S^2 \\
	& \hspace{10cm} + \mu_1 (2S^2 + 8A\varepsilon^2 + o(\varepsilon^2)) S^4\\
	& = \mu_2 S^2 (4S^4 + 8A S^2 \varepsilon^2 + o(\varepsilon^2))
	  + 2(128 D S^2 \varepsilon^2 + o(\varepsilon^2)) + \mu_1 (2S^2 + 8A\varepsilon^2 + o(\varepsilon^2)) S^4 \\
	& = 2(\mu_1 + 2\mu_2)S^6 + 8\{A(\mu_1 + \mu_2) + 32D\} S^4 \varepsilon^2 + o(\varepsilon^2).
\end{align*}
Thus, we obtain
\begin{equation}
\label{3.8}
\begin{aligned}
& \frac{\mu_2\|U_1\|^4_4 \|W_\varepsilon\|^4_{\mathcal{D}^{1,2}} + 2B(W_\varepsilon,U_1)\|W_\varepsilon\|_{\mathcal{D}^{1,2}}^2 \|U_1\|_{\mathcal{D}^{1,2}}^2 + \mu_1\|W_\varepsilon\|^4_4 \|U_1\|^4_{\mathcal{D}^{1,2}}}{2\mu_1 \mu_2 S^4} \\
&= \frac{\mu_1 + 2\mu_2}{\mu_1\mu_2} S^2 +\frac{4\{A(\mu_1 + \mu_2) + 32D\}}{\mu_1\mu_2}\varepsilon^2 + o(\varepsilon^2).
\end{aligned}
\end{equation}
Using Proposition \ref{prop 2.5}, \eqref{3.8}, and \eqref{3.7}, we obtain
\begin{align*}
	c_K & \le J(W_\varepsilon,U_1)\\
	& = \frac{1}{4}\frac{\mu_2\|U_1\|^4_4 \|W_\varepsilon\|^4_{\mathcal{D}^{1,2}} 
	  + 2B(W_\varepsilon,U_1)\|W_\varepsilon\|_{\mathcal{D}^{1,2}}^2 \|U_1\|_{\mathcal{D}^{1,2}}^2 
	  + \mu_1\|W_\varepsilon\|^4_4 \|U_1\|^4_{\mathcal{D}^{1,2}}}
	  {\mu_1 \mu_2 \|W_\varepsilon\|_4^4 \|U_1\|_4^4 - B(W_\varepsilon, U_1)^2}\\
	& = \frac{1}{4} \left(\frac{\mu_1 + 2\mu_2}{\mu_1\mu_2} S^2 
	   + \frac{4\{A(\mu_1 + \mu_2) + 32D\}}{\mu_1\mu_2}\varepsilon^2 + o(\varepsilon^2)\right) 
	  \left(1 - \frac{4A}{S^2} \varepsilon^2 + o(\varepsilon^2)\right)\\
	& = \frac{\mu_1 + 2\mu_2}{4\mu_1\mu_2}S^2 + \frac{1}{\mu_1\mu_2}[\{A(\mu_1 + \mu_2) + 32D\} 
	  - A(\mu_1 + 2\mu_2)]\varepsilon^2 + o(\varepsilon^2)\\
	& = \frac{1}{4}\left(\frac{1}{\mu_1} + \frac{2}{\mu_2}\right)S^2 + \frac{32D - \mu_2 A}{\mu_1\mu_2}\varepsilon^2 + o(\varepsilon^2).
\end{align*}
By Remark \ref{rmk 3.4} (i) and \eqref{1.4}, we deduce 
\begin{equation*}
	32D - \mu_2 A = 32\int_{\R^4} \frac{\beta(x)}{|x|^4} \, dx - 32  \pi^2 \mu_2 
	= 32 \left(\int_{\R^4} \frac{\beta(x)}{|x|^4} \, dx - \pi^2 \mu_2 \right) <0.
\end{equation*}
Thus, we conclude $c_K < \frac{1}{4}\left(\frac{1}{\mu_1} + \frac{2}{\mu_2}\right)S^2$. 
\end{proof}

\section{Existence of Kelvin--invariant positive solutions}
In this section, we prove Theorems \ref{thm 1.1}, \ref{thm 1.2}, and \ref{thm 1.4}. 
To this end, we classify Palais--Smale sequences for $I$ by applying the abstract profile decomposition theorem established in \cite{Sato--Suzuki}. 

For $y \in \R^4$ and $\lambda \in \R$, we define translation operator $\tau_y$ and dilation operator $\delta_\lambda$ on $\mathcal{D}^{1,2}(\R^4)$ as follows:
\begin{equation*}
(\tau_y u)(x) = u(x-y), \quad (\delta_\lambda u)(x) = 2^{\lambda} u(2^\lambda x).
\end{equation*}
Furthermore, the translation--dilation operator ${T}_{y,\lambda}$ is defined by ${T}_{y,\lambda} := \tau_y \circ \delta_\lambda$, that is,
\begin{equation*}
	({T}_{y,\lambda}u)(x) = 2^{\lambda} u(2^{\lambda}(x-y)).
\end{equation*}
For the operator $T$ on $\mathcal{D}^{1,2}(\R^4)$, we define the operator $\bm{T}$ on $\mathcal{D}$ by
\begin{equation*}
\bm{T}(u_1,u_2) = (T u_1, T u_2).
\end{equation*}
We set 
\begin{equation*}
G = \{\bm{T}_{y,\lambda} \mid y \in \R^4, \lambda \in \R\}.
\end{equation*}
For $y \in \R^4$ and $\lambda \in \R$, $\bm{T}_{y,\lambda}$ is a unitary operator on $\mathcal{D}$ and an isometry on $L^4(\R^4)^2$. 
We recall the following characterization of weak convergence to zero for translation--dilation operators.

\begin{lem}[{\cite[Lemma 2.2]{Sato--Suzuki}}]
\label{lem a}
For sequences $(y_k) \subset \mathbb{R}^4$ and $(\lambda_k) \subset \R$, the following are equivalent. 
\begin{itemize}
\item[(i)] ${T}_{y_k,\lambda_k} \rightharpoonup 0$ weakly in the weak operator topology.
\item[(ii)] $|y_k| + |\lambda_k| \to \infty.$
\end{itemize}
\end{lem}

We introduce the auxiliary functional $I_a$. 
For $a \in \R$ and $\bm{u}\in\mathcal{D}$, we define $I_a$ by
\begin{align*}
	I_a(\bm{u}) = \frac{1}{2} \int_{\R^4} (|\nabla u_1|^2 + |\nabla u_2|^2) \, dx 
	- \frac{\mu_1}{4} \int_{\R^4} u_1^4 \, dx - \frac{\mu_2}{4} \int_{\mathbb{R}^4} u_2^4 \, dx + \frac{a}{2} \int_{\R^4} u_1^2 u_2^2 \, dx.
\end{align*}
We also define
\begin{align*}
\mathcal{K}_K &= \{\bm{u} \in \mathcal{D}_K \mid I^\prime(\bm{u})=0\},\\
\mathcal{K}_a &= \{\bm{u} \in \mathcal{D} \mid I_a^\prime(\bm{u})=0\}.
\end{align*}

The following lemma provides a lower bound for the critical levels of the functionals $I$ and $I_a$.

\begin{lem}
\label{lem 4.1}
Set $\mathcal{I} = I$ or $I_a$ for $a \ge 0$. 
Let $\bm{u}$ be a critical point of $\mathcal{I}$. Then, the following hold.
\begin{itemize}
\item[(i)] If $u_i \not= 0$ for some $i = 1,2$, then we have $\displaystyle \mathcal{I}(\bm{u}) \ge \frac{1}{4\mu_i}S^2$. 
\item[(ii)] If $u_i \not= 0$ for both $i = 1,2$, then we have $ \displaystyle \mathcal{I}(\bm{u}) \ge \frac{1}{4}\left(\frac{1}{\mu_1} + \frac{1}{\mu_2}\right) S^2 $.
\end{itemize}
\end{lem}

\begin{proof}
Assume $u_i \not= 0$ for some $i = 1,2$. 
Testing the corresponding equation with $u_i$, we obtain
\begin{equation}
\label{4.1}
\|u_i\|_{\mathcal{D}^{1,2}}^2 \le \mu_i \|u_{i}\|_4^4.
\end{equation}
By the definition of the best Sobolev constant $S$ and \eqref{4.1}, 
\begin{equation}
\label{4.2}
S^2 \le \frac{\|u_i\|^4_{\mathcal{D}^{1,2}}}{\|u_i\|^4_4}\le \mu_i \|u_i\|^2_{\mathcal{D}^{1,2}}.
\end{equation}
On the other hand, since $\bm{u}$ is a critical point of $\mathcal{I}$, we have
\begin{equation}
\label{4.3}
	\mathcal{I}(\bm{u}) 
	= \mathcal{I}(\bm{u})-\frac{1}{4}\mathcal{I}'(\bm{u})\bm{u} 
	= \frac{1}{4}(\|u_1\|_{\mathcal{D}^{1,2}}^2 + \|u_2\|^2_{\mathcal{D}^{1,2}}).
\end{equation}
Combining \eqref{4.2} with \eqref{4.3}, we obtain (i). 
If $u_i \not= 0$ for both $i = 1,2$ \eqref{4.2} holds for each $i$. 
Thus, (ii) follows from \eqref{4.3}. 
\end{proof}

The following proposition provides a profile decomposition for Palais--Smale sequences of $I$.

\begin{prop}
\label{prop 4.2}
Suppose ($\beta$1)--($\beta$2) and let $(\bm{u}_k) \subset \mathcal{D}_K$ be a (PS)$_d$ sequence for $I$. Then, up to a subsequence, there exist $n_0 \in \N \cup \{0\}$, $\bm{w}^{(0)} \in \mathcal{K}_K$, and $\bm{g}_k^{(0)} = \bm{id}$, and for each $n \in \{1,\dots, n_0\}$, there exists $\beta^{(n)} \in \R$, $\bm{w}^{(n)} \in \mathcal{K}_{\beta^{(n)}} \setminus\{\bm{0}\}$, $(\bm{g}_k^{(n)}) \subset G$ such that
\begin{align}
\label{4.4}
&{\bm{g}_k^{(n)}}^* \bm{u}_k \rightharpoonup \bm{w}^{(n)}\quad \text{ weakly in } \mathcal{D} \text{ for all } n \in \{0,1,\dots,n_0\},\\
\label{4.5}
&{\bm{g}_k^{(n)}}^* \bm{g}_k^{(m)} \rightharpoonup \bm{0} \quad \text{ weakly for all } n,m \in \{0,1,\dots,n_0\} \text{ with } n \not= m,\\
\label{4.6}
&{\bm{g}_k^{(n)}}^* \bm{K} \bm{g}_k^{(m)} \rightharpoonup \bm{0} \quad \text{ weakly for all } n,m \in \{1,\dots,n_0\},\\
\label{4.7}
&\left\|\bm{u}_k - \bm{w}^{(0)} -\sum_{n=1}^{n_0} (\bm{g}_k^{(n)} \bm{w}^{(n)} + \bm{K} \bm{g}_k^{(n)} \bm{w}^{(n)})\right\|_{\mathcal{D}} \to 0,\\
\label{4.8}
&d = \lim_{k \to \infty} I(\bm{u}_k) = I(\bm{w}^{(0)}) + 2\sum_{n=1}^{n_0} I_{\beta^{(n)}}(\bm{w}^{(n)}).
\end{align}
\end{prop}

\begin{proof}
Let $(\bm{u}_k) \subset \mathcal{D}_K$ be a (PS)$_d$ sequence for $I$. Then, we have
\begin{equation*}
\|\bm{u}_k\|_{\mathcal{D}}^2 = 4I(\bm{u}_k) - I^\prime(\bm{u}_k)\bm{u}_k \le 4d + \|I^\prime(\bm{u}_k)\|_{\mathcal{D}^*} \|\bm{u}_k\|_{\mathcal{D}} + o(1).
\end{equation*}
Thus, $(\bm{u}_k)$ is bounded in $\mathcal{D}$. 
We use the profile decomposition theorem \cite[Theorem 1.8]{Sato--Suzuki}.
Since the operators in $G$ act diagonally on the two components,
the assumptions of this theorem are satisfied with $H=\mathcal{D}$, $D=G$, and $K=\bm{K}$ by \cite[Lemmas 2.2 and 2.3]{Sato--Suzuki}.
Therefore, by \cite[Theorem 1.8]{Sato--Suzuki},
up to a subsequence, there exist $n_0 \in \N \cup \{0, \infty\}$, $\bm{w}^{(0)} \in \mathcal{D}_K$, $\bm{g}_k^{(0)} = \bm{id}$, 
and for each $n \in \{1,\dots, n_0\}$, there exists $\bm{w}^{(n)} \in \mathcal{D} \setminus\{\bm{0}\}$ and $(\bm{g}_k^{(n)}) \subset G$ such that
\begin{align}
\label{4.9}
&{\bm{g}_k^{(n)}}^* \bm{u}_k \rightharpoonup \bm{w}^{(n)}\quad \text{ weakly in } \mathcal{D} \text{ for all } n \in \{0,1,\dots,n_0\},\\
\label{4.10}
&{\bm{g}_k^{(n)}}^* \bm{g}_k^{(m)} \rightharpoonup \bm{0} \quad \text{ weakly for all } n,m \in \{0,1,\dots,n_0\} \text{ with } n \not= m,\\
\label{4.11}
&{\bm{g}_k^{(n)}}^* \bm{K} \bm{g}_k^{(m)} \rightharpoonup \bm{0} \quad \text{ weakly for all } n,m \in \{1,\dots,n_0\},\\
\label{4.12}
&\|\bm{w}^{(0)}\|_{\mathcal{D}}^2 + 2 \sum_{n=1}^{n_0} \|\bm{w}^{(n)}\|_{\mathcal{D}}^2 \le \limsup_{k \to \infty} \|\bm{u}_k\|_{\mathcal{D}}^2.
\end{align}
Moreover, if $n_0 < \infty$, then
\begin{equation}
\label{4.13}
\bm{r}_k := \bm{u}_k - \bm{w}^{(0)} - \sum_{n=1}^{n_0} (\bm{g}_k^{(n)}\bm{w}^{(n)} + \bm{K}\bm{g}_k^{(n)}\bm{w}^{(n)}) \to \bm{0} \quad \text{in } L^4(\R^4)^2.
\end{equation}
Here, the $L^4$-convergence in \eqref{4.13}
follows from the $D$-weak convergence of the remainder in \cite[Theorem 1.8]{Sato--Suzuki} and \cite[Lemma 4.2]{Sato--Suzuki}.
We divide the proof into several steps.

\medskip

\noindent
\textbf{Step 1.} \textit{We have $n_0 < \infty$. Moreover, $\bm{w}^{(0)} \in \mathcal{K}_K$ and for all $n \in \{1,\dots,n_0\}$, there exists $\beta^{(n)} \in \R$ such that $\bm{w}^{(n)} \in \mathcal{K}_{\beta^{(n)}} \setminus \{\bm{0}\}$.}

\begin{proof}[Proof of Step 1.]
By the definition of $G$, for each $n \in \{1, \dots, n_0\}$, there exist $(y_k^{(n)}) \subset \R^4$ and 
$(\lambda_k^{(n)}) \subset \R$ such that $\bm{g}_k^{(n)} = \bm{T}_{y_k^{(n)},\lambda_k^{(n)}}$. 
Let $\bm{\varphi} \in \mathcal{D}$. 
Since $(\bm{u}_k)$ is a (PS)$_d$ sequence, for each $n \in \{0,1,\dots,n_0\}$, we have
\begin{align*}
	o(1) & = I^\prime(\bm{u}_k)(\bm{g}_k^{(n)} \bm{\varphi})\\
	& =  (\bm{u}_k, \bm{g}_k^{(n)} \bm{\varphi})_{\mathcal{D}} - \mu_1 \int_{\R^4} u_{1,k}^3 (g_k^{(n)} \varphi_1) \, dx 
	  - \mu_{2} \int_{\R^4} u_{2,k}^3 (g_k^{(n)} \varphi_2) \, dx \\
	& \hspace{2cm} + \int_{\R^4} \beta(x) u_{1,k} u_{2,k}^2 (g_k^{(n)} \varphi_1) \, dx 
	  + \int_{\R^4} \beta(x) u_{1,k}^2 u_{2,k} (g_k^{(n)} \varphi_2) \, dx \\
	& = ({\bm{g}_k^{(n)}}^* \bm{u}_k, \bm{\varphi})_{\mathcal{D}} 
	  - \mu_1 \int_{\R^4} ({g_k^{(n)}}^* u_{1,k})^3 \varphi_1 \, dx 
	  - \mu_2 \int_{\R^4} ({g_k^{(n)}}^* u_{2,k})^3 \varphi_2 \, dx  \\
	& \hspace{2cm} + \int_{\R^4} \beta(2^{-\lambda_k^{(n)}} x +y_k^{(n)}) ({g_k^{(n)}}^* u_{1,k}) ({g_k^{(n)}}^* u_{2,k})^2 \varphi_1 \, dx \\
	& \hspace{2cm} + \int_{\R^4} \beta(2^{-\lambda_k^{(n)}} x +y_k^{(n)}) ({g_k^{(n)}}^* u_{1,k})^2 ({g_k^{(n)}}^* u_{2,k}) \varphi_2 \, dx.
\end{align*}
By $\bm{g}_k^{(0)} = \bm{id} $ and ${\bm{g}_k^{(0)}}^* \bm{u}_k  \rightharpoonup \bm{w}^{(0)}$ weakly in $\mathcal{D}$, we obtain
\begin{align*}
	(\bm{w}^{(0)}, \bm{\varphi})_{\mathcal{D}}
	& - \mu_1 \int_{\R^4} (w_1^{(0)})^3 \varphi_1 \, dx 
	  - \mu_2 \int_{\R^4} (w_2^{(0)})^3 \varphi_2 \, dx \\
	& + \int_{\R^4} \beta(x) w_1^{(0)} (w_2^{(0)})^2 \varphi_1 \, dx 
	  + \int_{\R^4} \beta(x) (w_1^{(0)})^2 w_2^{(0)}\varphi_2 \, dx = 0.
\end{align*}
Hence, we have $\bm{w}^{(0)} \in \mathcal{K}_K$. 
Next, let $n \in \{1,\dots,n_0\}$. 
Since $\bm{g}_k^{(n)} \rightharpoonup \bm{0}$ weakly by \eqref{4.10}, 
Lemma \ref{lem a} implies that $|y_k^{(n)}| + |\lambda_k^{(n)}| \to \infty$. 
Passing to a subsequence, we may assume that either $|y_k^{(n)}| \to \infty$, 
or $|\lambda_k^{(n)}| \to \infty$ and $y_k^{(n)} \to y^{(n)}$ for some $y^{(n)} \in \R^4$. 
If $\lambda_k^{(n)} \to -\infty$ or $|y_k^{(n)}| \to \infty$, then $\beta(2^{-\lambda_k^{(n)}} x +y_k^{(n)}) \to \beta(\infty)= \beta(0)$ a.e. $x \in \R^4$. 
By \eqref{4.9}, we obtain
\begin{align*}
	(\bm{w}^{(n)}, \bm{\varphi})_{\mathcal{D}}
	& - \mu_1 \int_{\R^4} (w_1^{(n)})^3 \varphi_1 \, dx 
	  - \mu_2 \int_{\R^4} (w_2^{(n)})^3 \varphi_2 \, dx \\
	& + \int_{\R^4} \beta(0) w_1^{(n)} (w_2^{(n)})^2 \varphi_1 \, dx 
	  + \int_{\R^4} \beta(0) (w_1^{(n)})^2 w_2^{(n)}\varphi_2 \, dx = 0.
\end{align*}
Hence, we have $\bm{w}^{(n)} \in \mathcal{K}_{\beta^{(n)}} \setminus \{\bm{0}\}$ with $\beta^{(n)} = \beta(0)$. 
If $\lambda_k^{(n)} \to \infty$ and $y_k^{(n)} \to y^{(n)}$, then $\beta(2^{-\lambda_k^{(n)}} x +y_k^{(n)}) \to \beta(y^{(n)})$ a.e. $x \in \R^4$. By \eqref{4.9}, we have
\begin{align*}
	(\bm{w}^{(n)}, \bm{\varphi})_{\mathcal{D}}
	& - \mu_1 \int_{\R^4} (w_1^{(n)})^3 \varphi_1 \, dx 
	  - \mu_2 \int_{\R^4} (w_2^{(n)})^3 \varphi_2 \, dx \\
	& + \int_{\R^4} \beta(y^{(n)}) w_1^{(n)} (w_2^{(n)})^2 \varphi_1 \, dx 
	  + \int_{\R^4} \beta(y^{(n)}) (w_1^{(n)})^2 w_2^{(n)}\varphi_2 \, dx = 0.
\end{align*}
Hence, we have $\bm{w}^{(n)} \in \mathcal{K}_{\beta^{(n)}} \setminus \{\bm{0}\}$ with $\beta^{(n)} = \beta(y^{(n)})$. For each $n \in \{1,\dots,n_0\}$, by Lemma \ref{lem 4.1} and $ \bm{w}^{(n)} \not= \bm{0}$, there exists $i \in \{1,2\}$ such that
\begin{equation*}
	I_{\beta^{(n)}}(\bm{w}^{(n)}) 
	= \frac{1}{4}(\|w_1^{(n)} \|_{\mathcal{D}^{1,2}}^2 + \|w_2^{(n)}\|_{\mathcal{D}^{1,2}}^2) 
	\ge \frac{1}{4\mu_i}S^2.
\end{equation*}
Combining this with \eqref{4.12}, we conclude that $n_0 < \infty$.
\end{proof}

\noindent
\textbf{Step 2.} \textit{\eqref{4.7} holds.}
\begin{proof}[Proof of Step 2]
By $n_0 < \infty$, \eqref{4.13} holds. 
Since $(\bm{u}_k)$ is a (PS)$_d$ sequence for $I$, by \eqref{4.13}, we have
\begin{equation}
\label{4.14}
\begin{aligned}
	o(1) = I^\prime(\bm{u}_k)\bm{r}_k 
	& = (\bm{u}_k,\bm{r}_k)_{\mathcal{D}} - \mu_1 \int_{\R^4}  u_{1,k}^3 r_{1,k} \, dx -  \mu_2 \int_{\R^4}  u_{2,k}^3 r_{2,k} \, dx \\
	& \hspace{2cm} + \int_{\R^4} \beta(x) u_{1,k} u_{2,k}^2 r_{1,k} \, dx + \int_{\R^4} \beta(x) u_{1,k}^2 u_{2,k} r_{2,k} \, dx \\
	& = (\bm{u}_k,\bm{r}_k)_{\mathcal{D}} + o(1).
\end{aligned}
\end{equation}
From $\bm{w}^{(0)} \in \mathcal{K}_K$, as in \eqref{4.14}, we obtain
\begin{equation}
\label{4.15}
0 = I^\prime(\bm{w}^{(0)}) \bm{r}_k = (\bm{w}^{(0)}, \bm{r}_k)_{\mathcal{D}} + o(1).
\end{equation}
Moreover, for each $n \in \{1,\dots,n_0\}$, since $\bm{w}^{(n)} \in \mathcal{K}_{\beta^{(n)}}$ and $I_{\beta^{(n)}}$ 
is invariant  conformal transformations, as in \eqref{4.14}, we have
\begin{align}
\label{4.16}
&0 = I^\prime_{\beta^{(n)}}(\bm{w}^{(n)})({\bm{g}_k^{(n)}}^* \bm{r}_k) = I^\prime_{\beta^{(n)}}(\bm{g}_k^{(n)} \bm{w}^{(n)})\bm{r}_k = (\bm{g}_k^{(n)} \bm{w}^{(n)}, \bm{r}_k)_{\mathcal{D}} + o(1),\\
\label{4.17}
&0 = I^\prime_{\beta^{(n)}}(\bm{w}^{(n)})({\bm{g}_k^{(n)}}^* \bm{K} \bm{r}_k) = I^\prime_{\beta^{(n)}}(\bm{K}\bm{g}_k^{(n)} \bm{w}^{(n)})\bm{r}_k = (\bm{K}\bm{g}_k^{(n)} \bm{w}^{(n)},\bm{r}_k)_{\mathcal{D}} + o(1).
\end{align}
Combining \eqref{4.14}--\eqref{4.17}, we obtain
\begin{align*}
o(1) &= I^\prime(\bm{u}_k)\bm{r}_k - I^\prime(\bm{w}^{(0)}) \bm{r}_k - \sum_{n=1}^{n_0} \{I^\prime_{\beta^{(n)}}(\bm{g}_k^{(n)} \bm{w}^{(n)})\bm{r}_k + I^\prime_{\beta^{(n)}}(\bm{K}\bm{g}_k^{(n)} \bm{w}^{(n)})\bm{r}_k\}\\
&= (\bm{u}_k,\bm{r}_k)_{\mathcal{D}} - (\bm{w}^{(0)}, \bm{r}_k)_{\mathcal{D}} - \sum_{n=1}^{n_0} \{(\bm{g}_k^{(n)} \bm{w}^{(n)}, \bm{r}_k)_{\mathcal{D}} + (\bm{K}\bm{g}_k^{(n)} \bm{w}^{(n)},\bm{r}_k)_{\mathcal{D}}\} + o(1)\\
&= \|\bm{r}_k\|^2_{\mathcal{D}} + o(1).
\end{align*}
Therefore, \eqref{4.7} holds.
\end{proof}

\noindent
\textbf{Step 3.} \textit{\eqref{4.8} holds.}
\begin{proof}[Proof of Step 3]
First, by \eqref{4.7} and \eqref{4.9}--\eqref{4.11}, we have
\begin{equation}
\label{4.18}
\lim_{k\to\infty}\|\bm{u}_k\|^2_{\mathcal{D}} = \|\bm{w}^{(0)}\|_{\mathcal{D}}^2 + 2 \sum_{n=1}^{n_0} \|\bm{w}^{(n)}\|_{\mathcal{D}}^2.
\end{equation}
Next, for both $i = 1,2$, it follows that
\begin{equation}
\label{4.19}
\lim_{k\to\infty} \|u_{i,k}\|_4^4 =  \|w_i^{(0)}\|_4^4 + 2 \sum_{n=1}^{n_0} \|w_i^{(n)}\|_4^4.
\end{equation}
Indeed, using \eqref{4.9}--\eqref{4.11} and Brezis--Lieb lemma  \cite{Brezis--Lieb}, we deduce
\begin{align*}
\|r_{i,k}\|_4^4 &= \|u_{i,k}\|_4^4 - \|w_i^{(0)}\|_4^4 - \sum_{n=1}^{n_0} \|g_k^{(n)} w_i^{(n)}\|_4^4 - \sum_{n=1}^{n_0} \|K g_k^{(n)} w_i^{(n)}\|_4^4 + o(1)\\
&= \|u_{i,k}\|_4^4 - \|w_i^{(0)}\|_4^4 - 2 \sum_{n=1}^{n_0} \|w_i^{(n)}\|_4^4 + o(1).
\end{align*}
Combining this with \eqref{4.13}, we obtain \eqref{4.19}. 
Furthermore, we have
\begin{equation}
\label{4.20}
\lim_{k\to\infty} B(u_{1,k},u_{2,k}) = B(w_1^{(0)},w_2^{(0)}) + 2\sum_{n=1}^{n_0} \beta^{(n)} \|w_1^{(n)} w_2^{(n)}\|_2^2.
\end{equation}
In fact, using \eqref{4.9}--\eqref{4.11} and \eqref{4.13}, we have
\begin{equation*}
\label{A}
\begin{aligned}
	& \quad B(u_{1,k},u_{2,k}) \\
	& = B\left(w_1^{(0)} + \sum_{n=1}^{n_0}\left(g_k^{(n)} w_1^{(n)} 
	  + Kg_k^{(n)} w_1^{(n)}\right) + r_{1,k}, w_2^{(0)} 
	  + \sum_{n=1}^{n_0}\left(g_k^{(n)} w_2^{(n)} 
	  + Kg_k^{(n)} w_2^{(n)}\right) + r_{2,k} \right)\\
	& = B(w_1^{(0)}, w_2^{(0)}) + \sum_{n=1}^{n_0} \left\{B(g_k^{(n)} w_1^{(n)},g_k^{(n)} w_2^{(n)}) 
	  + B(Kg_k^{(n)} w_1^{(n)}, Kg_k^{(n)} w_2^{(n)})\right\} + o(1).
\end{aligned}
\end{equation*}
By a change of variables and the definition of $\beta^{(n)}$, we obtain
\begin{equation*}
\begin{aligned}
	B(Kg_k^{(n)} w_1^{(n)}, Kg_k^{(n)} w_2^{(n)}) 
	& = B(g_k^{(n)} w_1^{(n)}, g_k^{(n)} w_2^{(n)}) \\
	& = \int_{\R^4} \beta(2^{-\lambda_k^{(n)}}x + y_k^{(n)})  ({w_1^{(n)}})^2 ({w_2^{(n)}})^2 \, dx \\
	& = \beta^{(n)} \|{w_1^{(n)}} {w_2^{(n)}}\|_2^2 + o(1)
\end{aligned}
\end{equation*}
Combining this with \eqref{4.13}, we have \eqref{4.20}. By \eqref{4.18}--\eqref{4.20}, \eqref{4.8} holds.
\end{proof}
Steps 1--3 complete the proof.

\end{proof}

\begin{proof}[Proof of Theorem \ref{thm 1.1}]
By Lemma \ref{lem 2.3.5}, there exists a (PS)$_{c_K}$ sequence $(\bm{u}_k)\subset\mathcal{N}_K$ for $I$.
By Proposition \ref{prop 4.2}, 
up to a subsequence, there exist $n_0 \in \N \cup \{0\}$, $\bm{w}^{(0)}\in \mathcal{K}_K$, and 
$\bm{g}_k^{(0)} = \bm{id}$, and for each $n \in \{1,\dots, n_0\}$, there exist $\beta^{(n)} \in \R$, 
$\bm{w}^{(n)} \in \mathcal{K}_{\beta^{(n)}} \setminus\{\bm{0}\}$, 
$(\bm{g}_k^{(n)}) \subset G$ such that
\begin{align}
	& {\bm{g}_k^{(n)}}^* \bm{u}_k \rightharpoonup \bm{w}^{(n)}\quad \text{ weakly in } \mathcal{D} \text{ for all } n \in \{0,1,\dots,n_0\},\\
	& {\bm{g}_k^{(n)}}^* \bm{g}_k^{(m)} \rightharpoonup \bm{0} \quad \text{ weakly for all } n,m \in \{0,1,\dots,n_0\} \text{ with } n \not= m,\\
	& {\bm{g}_k^{(n)}}^* \bm{K} \bm{g}_k^{(m)} \rightharpoonup \bm{0} \quad \text{ weakly for all } n,m \in \{1,\dots,n_0\},\\
	\label{4.26}&\left\|\bm{u}_k - \bm{w}^{(0)} -\sum_{n=1}^{n_0} (\bm{g}_k^{(n)} \bm{w}^{(n)} + \bm{K} \bm{g}_k^{(n)} \bm{w}^{(n)})\right\|_{\mathcal{D}} \to 0,\\
	\label{4.27} & c_K = \lim_{k \to \infty} I(\bm{u}_k) = I(\bm{w}^{(0)}) + 2\sum_{n=1}^{n_0} I_{\beta^{(n)}}(\bm{w}^{(n)}).
\end{align}
We first show that $w_1^{(0)}\not=0$ and $w_2^{(0)}\not=0$.
Suppose that $w_1^{(0)}=0$. 
Then, by Lemma \ref{lem 2.2} and \eqref{4.26}, at least one of $\bm{w}^{(1)},\ldots,\bm{w}^{(n_0)}$ has a nonzero first component, 
and at least one of $\bm{w}^{(0)},\bm{w}^{(1)},\ldots,\bm{w}^{(n_0)}$ has a nonzero second component.
Hence, by Lemma \ref{lem 4.1} and \eqref{4.27}, we have
\begin{equation*}
c_K  \ge \frac{1}{4}\left(\frac{2}{\mu_1}+\frac{1}{\mu_2}\right)S^2.
\end{equation*}
Similarly, if $w_2^{(0)}=0$, then 
\begin{equation*}
c_K  \ge \frac{1}{4}\left(\frac{1}{\mu_1}+\frac{2}{\mu_2}\right)S^2.
\end{equation*}
Since $\mu_1\leq \mu_2$, both inequalities contradict Proposition \ref{prop 3.1}. 
Hence, we have $w_1^{(0)}\not=0$ and $w_2^{(0)}\not=0$.
Suppose that $n_0\geq 1$. 
Since both components of $\bm{w}^{(0)}$ are nonzero, Lemma \ref{lem 4.1} yields
$I(\bm{w}^{(0)}) \geq \frac14\left(\frac1{\mu_1}+\frac1{\mu_2}\right)S^2$.
Moreover, since $\bm{w}^{(1)}\in\mathcal K_{\beta^{(1)}}\setminus\{\bm{0}\}$ and $\mu_1\leq \mu_2$, 
Lemma \ref{lem 4.1} (i) gives $I_{\beta^{(1)}}(\bm{w}^{(1)}) \ge \frac{1}{4\mu_2} S^2$. 
Therefore, by \eqref{4.27}, we have
\begin{equation*}
c_K = I(\bm{w}^{(0)}) + 2\sum_{n=1}^{n_0} I_{\beta^{(n)}}(\bm{w}^{(n)}) \ge \frac{1}{4}\left(\frac{1}{\mu_1}+\frac{3}{\mu_2}\right)S^2.
\end{equation*}
This inequality also contradicts Proposition \ref{prop 3.1}. 
Thus, $n_0=0$. By \eqref{4.26}, we see that $\bm{u}_k \to \bm{w}^{(0)}$ in $\mathcal{D}$. Therefore, $\bm{w}^{(0)} \in \mathcal{D}_K$ is a minimizer for $c_K$ and a Kelvin-invariant solution of \eqref{P}. Replacing $w_1^{(0)}, w_2^{(0)}$ by $|w_1^{(0)}|, |w_2^{(0)}|$ respectively, we may assume that $w_1^{(0)}, w_2^{(0)} \ge 0$. By the strong maximum principle, $w_1^{(0)}, w_2^{(0)}$ are positive. Therefore, the proof is complete.
\end{proof}

\begin{proof}[Proof of Theorem \ref{thm 1.2}]
It suffices to use Proposition \ref{prop 3.2} instead of Proposition \ref{prop 3.1} in the proof of Theorem \ref{thm 1.1}.
\end{proof}

\begin{proof}[Proof of Theorem \ref{thm 1.4}]
Since $\beta(x)$ is a radial function, we can see that $\mathcal{N}_K \cap \mathcal{D}_{\rm rad}$ 
is a natural constraint for $I$, by the same argument as Lemma \ref{lem 2.3}. 
Furthermore, since $(U_1,U_1)$ and $(U_\varepsilon + U_{\frac{1}{\varepsilon}},U_1)$ are radially symmetric, 
the same energy estimates for $c_{K,{\rm rad}}$ follow as in the proofs of Propositions \ref{prop 3.1} and \ref{prop 3.2}.
The rest of the proof is analogous to that of Theorem \ref{thm 1.1}.
\end{proof}

\section{Asymptotic behavior of minimizers}
\noindent
In this section, we consider the asymptotic behavior of the minimizers of $c_K$ and $c_{K,{\rm rad}}$ as $\|\beta\|_\infty\to 0$. 
Since the arguments for $c_{K,{\rm rad}}$ are essentially the same as those for $c_K$, we only give the proof for $c_K$. 
In what follows, we write $I=I_\beta$, $B = B_\beta$, $J=J_\beta$, $c_K=c_K(\beta)$, $\mathcal{N}_K=\mathcal{N}_{K}(\beta)$, 
$\mathcal{X}_K = \mathcal{X}_K(\beta)$ to highlight the dependence on $\beta$.

The following lemma is useful for studying the asymptotic behavior.

\begin{lem}
The following hold.
\label{lem 5.1}
\begin{itemize}
\item[(i)] If $0 \le \beta_1(x) \le \beta_2(x)$ for all $x \in \R^4$, then $c_K(\beta_1) \le c_K(\beta_2)$.
\item[(ii)] $\displaystyle c_K(0) = \frac{1}{4}\left(\frac{1}{\mu_1}+\frac{1}{\mu_2}\right)S^2$.
\item[(iii)] $c_K(\beta)\to c_K(0)$ as $\|\beta\|_\infty \to 0$.
\end{itemize}
\end{lem}

\begin{proof}
If $0 \le \beta_1(x) \le \beta_2(x)$ for all $x \in \R^4$, then we have $\mathcal{X}_K(\beta_2) \subset \mathcal{X}_K(\beta_1)$ and
\begin{equation*}
J_{\beta_1}(\bm{u}) \le J_{\beta_2}(\bm{u}) \quad \text{for all } \bm{u} \in \mathcal{X}_K(\beta_2) .
\end{equation*}
Hence, we obtain
\begin{equation*}
	c_K(\beta_1)
	=\inf_{\bm{u}\in \mathcal{X}_K(\beta_1)}J_{\beta_1}(\bm{u})
	\leq \inf_{\bm{u}\in \mathcal{X}_K(\beta_2)}J_{\beta_1}(\bm{u})
	\leq \inf_{\bm{u}\in \mathcal{X}_K(\beta_2)}J_{\beta_2}(\bm{u})
	=c_K(\beta_2).
\end{equation*}
Thus, (i) holds.  
Next, we prove (ii).
Note that $ \mathcal{X}_K(0) = \left\{\bm{u}\in \mathcal{D}_K\,\middle|\,u_1\not=0, u_2\not=0\right\}$.
Since the best constant in the Sobolev inequality is attained by the Kelvin-invariant Talenti function $U$, we have
\begin{equation*}
	S^2 = \inf_{u \in \mathcal{D}^{1,2}(\R^4) \setminus\{0\}} \frac{\|u\|_{\mathcal{D}^{1,2}}^4}{ \|u\|_4^4} 
	= \inf_{\substack{u \in \mathcal{D}^{1,2}(\R^4) \setminus\{0\},\\ Ku=u}} \frac{\|u\|_{\mathcal{D}^{1,2}}^4}{ \|u\|_4^4}.
\end{equation*}
By \eqref{2.3} with $\beta = 0$, we have
\begin{align*}
c_K(0) &= \frac{1}{4} \inf_{(u_1,u_2) \in \mathcal{X}_K(0)} \left( \frac{\|u_1\|_{\mathcal{D}^{1,2}}^4}{\mu_1 \|u_1\|_4^4} + \frac{\|u_2\|_{\mathcal{D}^{1,2}}^4}{\mu_2 \|u_2\|_4^4} \right)\\
&= \frac{1}{4} \left(\frac{1}{\mu_1}\inf_{u_1 \in \mathcal{D}^{1,2}(\R^4) \setminus\{0\}} \frac{\|u_1\|_{\mathcal{D}^{1,2}}^4}{ \|u_1\|_4^4} + \frac{1}{\mu_2} \inf_{u_2 \in \mathcal{D}^{1,2}(\R^4) \setminus\{0\}} \frac{\|u_2\|_{\mathcal{D}^{1,2}}^4}{ \|u_2\|_4^4}\right)\\
&= \frac{1}{4}\left(\frac{1}{\mu_1}+\frac{1}{\mu_2}\right)S^2.
\end{align*}
Therefore, we conclude (ii).
Finally, we prove (iii). 
Let $(\beta_k)$ satisfy $\|\beta_k\|_\infty \to 0$. 
Since $B_{\beta_k}(U,U)=o(1)$ as $k\to \infty$, we have $(U,U) \in \mathcal{X}_K(\beta_k)$ for large $k$.
Then, by \eqref{2.3}, we have
\begin{equation*}
	\limsup_{k\to\infty}c_K(\beta_k) \leq \limsup_{k\to\infty}J_{\beta_k}(U,U) = J_0(U,U)=c_K(0)
\end{equation*}
On the other hand, by Lemma \ref{lem 5.1} (i), we obtain
\begin{equation*}
	\liminf_{k\to\infty}c_K(\beta_k) \ge c_K(0).
\end{equation*}
Hence, (iii) holds.
\end{proof}

\begin{proof}[Proof of Theorem \ref{thm 1.6}]
We prove the assertion for positive minimizers of $c_K$. 
The proof for positive minimizers of $c_{K,{\rm rad}}$ is analogous and hence omitted.
Let $(\beta_k)$ be a sequence such that $\|\beta_k\|_\infty \to 0$ and let $\bm{w}_k $ be a minimizer for $c_K(\beta_k)$. 
Then $I_{\beta_k}(\bm{w}_k)=c_K(\beta_k)$ and $I_{\beta_k}'(\bm{w}_k)=0$.
By Lemma \ref{lem 5.1} (iii), we have
\begin{equation*}
\|\bm{w}_k\|_{\mathcal{D}}^2 = 4I_{\beta_k}(\bm{w}_k) - I_{\beta_k}^\prime(\bm{w}_k)\bm{w}_k = 4 c_K(\beta_k) \to 4 c_K(0).
\end{equation*}
Thus, $(\bm{w}_k)$ is bounded in $\mathcal{D}$. 
Hence, for $\bm{\varphi}= (\varphi_1,\varphi_2)\in \mathcal{D}$, we have
\begin{align*}
	|B_{\beta_k}(\bm{w}_k)| & \leq \|\beta_k\|_\infty\|w_{1,k}\|_{4}^2\|w_{2,k}\|_{4}^2=o(1),\\
	|B_{\beta_k}'(\bm{w}_k)\bm{\varphi}| 
	& = \left|2\int_{\mathbb{R}^4}\beta_k(x)\left(w_{1,k}w_{2,k}^2\varphi_1 + w_{1,k}^2 w_{2,k}\varphi_2\right) dx \right|\\
	& \leq 4S^{-2}\|\beta_k\|_\infty\|\bm{w}_k\|_{\mathcal{D}}^3\|\bm{\varphi}\|_{\mathcal{D}}  \\
	& = o(1)\|\bm{\varphi}\|_{\mathcal{D}}
\end{align*}
By Lemma \ref{lem 5.1} (iii), we deduce
\begin{align*}
I_0(\bm{w}_k) & = I_{\beta_k}(\bm{w}_k) + \frac{1}{2}B_{\beta_k}(\bm{w}_k) = c_K(\beta_k) + o(1)\to c_K(0),\\
\|I_0^\prime(\bm{w}_k)\|_{\mathcal{D}^*}
& \leq \|I^\prime_{\beta_k}(\bm{w}_k)\|_{\mathcal{D}^*} + \frac{1}{2}\|B_{\beta_k}'(\bm{w}_k)\|_{\mathcal{D}^*} = o(1).
\end{align*}
Hence, $(\bm{w}_k)$ is a (PS)$_{c_K(0)}$ sequence for $I_0$. 
Then, by Proposition \ref{prop 4.2}, up to a subsequence, there exist $n_0 \in \N \cup \{0\}$, $\bm{w}^{(0)}\in \mathcal{K}_K$, and 
$\bm{g}_k^{(0)} = \bm{id}$, and for each $n \in \{1,\dots, n_0\}$, there exists $\bm{w}^{(n)} \in \mathcal{K}_0 \setminus\{\bm{0}\}$, 
$(\bm{g}_k^{(n)}) \subset G$ such that
\begin{align}
	& {\bm{g}_k^{(n)}}^* \bm{w}_k \rightharpoonup \bm{w}^{(n)}\quad \text{ weakly in } \mathcal{D} \text{ for all } n \in \{0,1,\dots,n_0\},\\
	& {\bm{g}_k^{(n)}}^* \bm{g}_k^{(m)} \rightharpoonup \bm{0} \quad \text{ weakly for all } n,m \in \{0,1,\dots,n_0\} \text{ with } n \not= m,\\
	& {\bm{g}_k^{(n)}}^* \bm{K} \bm{g}_k^{(m)} \rightharpoonup \bm{0} \quad \text{ weakly for all } n,m \in \{1,\dots,n_0\},\\
	\label{5.4}&\left\|\bm{w}_k - \bm{w}^{(0)} -\sum_{n=1}^{n_0} (\bm{g}_k^{(n)} \bm{w}^{(n)} + \bm{K} \bm{g}_k^{(n)} \bm{w}^{(n)})\right\|_{\mathcal{D}} \to 0,\\
	\label{5.5} & c_K(0) = \lim_{k \to \infty} I_0(\bm{w}_k) = I_0(\bm{w}^{(0)}) + 2\sum_{n=1}^{n_0} I_0(\bm{w}^{(n)}).
\end{align}
We first show that $w_1^{(0)}\not=0$ and $w_2^{(0)}\not=0$.
Suppose that $w_1^{(0)}=0$. 
Then, by Lemma \ref{lem 2.2} and \eqref{5.4}, at least one of $\bm{w}^{(1)},\ldots,\bm{w}^{(n_0)}$ has a nonzero first component, 
and at least one of $\bm{w}^{(0)},\bm{w}^{(1)},\ldots,\bm{w}^{(n_0)}$ has a nonzero second component.
Hence, by Lemma \ref{lem 4.1} and \eqref{5.5}, we have
\begin{equation*}
c_K(0)  \ge \frac{1}{4}\left(\frac{2}{\mu_1}+\frac{1}{\mu_2}\right)S^2.
\end{equation*}
Similarly, if $w_2^{(0)}=0$, then 
\begin{equation*}
c_K(0)  \ge \frac{1}{4}\left(\frac{1}{\mu_1}+\frac{2}{\mu_2}\right)S^2.
\end{equation*}
Both inequalities contradict Lemma \ref{lem 5.1} (ii). 
Hence, we have $w_1^{(0)}\not=0$ and $w_2^{(0)}\not=0$.
Suppose that $n_0\geq 1$. 
Since both components of $\bm{w}^{(0)}$ are nonzero, Lemma \ref{lem 4.1} yields
$I_0(\bm{w}^{(0)}) \geq \frac14\left(\frac1{\mu_1}+\frac1{\mu_2}\right)S^2$.
Moreover, since $\bm{w}^{(1)}\in\mathcal K_0\setminus\{\bm{0}\}$, 
Lemma \ref{lem 4.1} gives $I_0(\bm{w}^{(1)})>0$. 
Therefore, by \eqref{5.5}, we have
\begin{equation*}
c_K(0) = I_0(\bm{w}^{(0)}) + 2\sum_{n=1}^{n_0} I_0(\bm{w}^{(n)}) > \frac{1}{4}\left(\frac{1}{\mu_1}+\frac{1}{\mu_2}\right)S^2=c_K(0),
\end{equation*}
which is a contradiction. 
Thus, $n_0=0$.
By \eqref{5.4}, we see that $\bm{w}_k \to \bm{w}^{(0)}$ in $\mathcal{D}$. 
Hence, $\bm{w}^{(0)}$ is a critical point for $I_0$ satisfying $I_0(\bm{w}^{(0)}) = c_K(0)$ and $\bm{w}^{(0)} \in \mathcal{D}_K$. 
Hence, $w^{(0)}_i$ is a Kelvin invariant solution to \eqref{1.6} with $\mu = \mu_i$.
Since $\bm{w}_k$ is positive, we have $w_i^{(0)}\geq 0$ for $i=1,2$.
As both components are nontrivial, the strong maximum principle implies that they are positive.
By the uniqueness of positive Kelvin-invariant solutions of \eqref{1.6},
we have $w^{(0)}_i  = \mu_{i}^{-\frac{1}{2}} U$. 
Thus, we obtain 
\begin{equation}
\label{5.6}
\bm{w}_k \to (\mu_1^{-\frac{1}{2}}U, \mu_2^{-\frac{1}{2}}U) \quad \text{in }\mathcal{D}_K.
\end{equation}
Since the above argument applies to every subsequence, \eqref{5.6} holds with the original sequence. Therefore, the proof is complete.
\end{proof}

\subsection*{Acknowledgements}
This work was supported by JSPS KAKENHI Grant Numbers JP20K03691.

\subsection*{Data Availability}
No data was used for the research described in the article.

\end{document}